\documentclass[11pt,a4paper,reqno]{amsart}

\usepackage[T1]{fontenc}
\usepackage[utf8]{inputenc}
\usepackage{lmodern}
\usepackage[margin=30mm]{geometry}
\usepackage{microtype}
\usepackage{amsmath,amssymb,amsthm,mathtools,mathrsfs}
\usepackage{booktabs}
\usepackage{xcolor}
\usepackage{hyperref}
\usepackage{bookmark}
\hypersetup{%
  colorlinks=true, linkcolor=blue!60!black, citecolor=blue!60!black, urlcolor=blue!60!black,
  pdftitle={Stability and bifurcation of shrinking spheres for the fractional mean curvature flow},
  pdfauthor={Annalisa Cesaroni, Matteo Novaga},
  pdfsubject={Stability above the threshold, instability below, and near-threshold shrinkers},
  pdfkeywords={fractional mean curvature flow, shrinking spheres, nonlinear stability, bifurcation, self-similar shrinkers}
}
\numberwithin{equation}{section}
\allowdisplaybreaks[2]
\newtheorem{theorem}{Theorem}[section]
\newtheorem{proposition}[theorem]{Proposition}
\newtheorem{lemma}[theorem]{Lemma}
\newtheorem{corollary}[theorem]{Corollary}
\newtheorem{maintheorem}{Theorem}

\theoremstyle{definition}
\newtheorem{definition}[theorem]{Definition}

\theoremstyle{remark}

\newcommand{\R}{\mathbb R}
\newcommand{\Sph}{{\mathbb S^{n-1}}}
\newcommand{\dd}{\,\mathrm d}
\newcommand{\pv}{\operatorname{PV}}
\newcommand{\Hs}{H_s}

\newcommand{\Is}{\mathcal I_s}

\title[Stability and bifurcation of shrinking spheres]{Stability and bifurcation of shrinking spheres \\ for the fractional mean curvature flow}

\author{Annalisa Cesaroni}
\address{Dipartimento di Matematica ``Tullio Levi-Civita'', Università di Padova, Via Trieste 63, 35131 Padova, Italy}
\email{annalisa.cesaroni@unipd.it}
\author{Matteo Novaga}
\address{Dipartimento di Matematica, Università di Pisa, Largo Bruno Pontecorvo 5, 56127 Pisa, Italy}
\email{matteo.novaga@unipi.it}

\date{}


\begin{document}
\begin{abstract}
We study the fractional mean curvature flow near shrinking spheres. The first shape eigenvalue changes sign at
\(s_*(n)=(\sqrt{n^2+6n+1}-n-1)/2\). Above this threshold, sufficiently small perturbations of a sphere become extinct at one point, and their extinction rescalings converge exponentially to the round shrinker. Below the threshold, the sphere is nonlinearly unstable modulo translations and dilations. We also construct smooth, strictly convex, nonspherical homothetic shrinkers bifurcating from the spherical family near the threshold. In the plane, the noncircular branch bifurcates on the side \(s>s_*(2)\). \end{abstract}
\maketitle
\tableofcontents

\section{Introduction and main results}
\label{sec:problem}

 The sphere is the basic self-similar solution of fractional mean curvature flow. In this paper we analyse the dynamical stability of this round shrinker under the unconstrained flow. Moreover, we discuss what happens at the first loss of stability, and whether nonspherical convex shrinkers bifurcate from the spherical family. 
 
 These questions differ in an essential way from stability for the volume-preserving fractional flow: extinction rescaling of the unconstrained evolution produces a dilation term, and the competition between that term and the nonlocal Jacobi operator leads to a dimension-dependent spectral threshold. 
 
Evolution and first-variation formulas for fractional mean curvature flow were developed in \cite{SV}, and classical short-time existence was proved in \cite{JLMshort}; the level-set formulation and viscosity solution approach go back to  \cite{Imbert}.   Broader weak and variational theories for nonlocal curvature motions were developed in  \cite{CMP}. We recall here also   \cite{CNR} where  anisotropic fractional mean curvature flows are studied through threshold dynamics, including consistency, preservation of convexity, and uniqueness for bounded convex evolutions. Stability of nonlocal geometric flows under convergence of the curvatures, including the limits of fractional mean curvature flows as $s\to0^+$ and $s\to1^-$, was studied in \cite{CDNP}. For graphical geometries, the authors in \cite{CNgraphs} obtained regularity and long-time asymptotics for entire Lipschitz graphs, while the paper \cite{MWgraphs} developed a quasilinear theory for periodic graphs.

The unconstrained flow can develop singularities other than spherical extinction: in \cite{CSVneck}   embedded initial hypersurfaces developing neckpinch singularities before shrinking to a point are constructed. 

A different line of work concerns the volume-preserving velocity $-H_s+\overline H_s$. For convex sets, a priori estimates and asymptotic roundness were established in \cite{CSVvolume}; near a ball, long-time existence and asymptotic stability were proved in \cite{CNball}, by a  quantitative Alexandrov type estimate for
nearly spherical sets.  Finally, \cite{JLMconvex} proves global regularity and convergence to a ball for convex initial data under the volume-preserving flow.

These results do not imply the analysis provided  in this paper: after extinction rescaling the unconstrained flow has velocity $-H_s+(1+s)^{-1}x\cdot\nu$, and the dilation term changes the shape spectrum. The sign change of its degree-two eigenvalue determines the threshold $s=s_*(n)$ defined below.

In the classical curve shortening flow, the circle is the only embedded closed homothetically shrinking curve, by the classification of Abresch and Langer \cite{AbreschLanger}, and convex curves shrink to round points by the theorem of Gage and Hamilton \cite{GageHamilton}; for hypersurfaces, Huisken \cite{Huisken} proved that the only compact mean-convex self-shrinkers are spheres. In this paper we  show that this rigidity fails for the fractional flow near $s=s_*(n)$, already among connected, strictly convex shrinkers.

We recall that self-similar solutions of the fractional flow were previously constructed in \cite{CNshrink}: their boundaries consist of a prescribed number of concentric spheres, and all the nonspherical examples in that family are dynamically unstable. The branches constructed here are of a different type: their boundary is connected and is a small, strictly convex deformation of one sphere.

The shrinker equation must be distinguished from the equation of constant nonlocal mean curvature. The nonlocal Alexandrov theorems proved independently in \cite{CFSMW} and in  \cite{CFMN} show that a bounded smooth set with constant fractional mean curvature is a ball; the periodic Delaunay-type sets of \cite{CFSMW,CFWdelaunay} and the near-sphere lattices of \cite{CFWlattices} solve that  equation. On the nonspherical shrinkers constructed here the fractional curvature is not constant: it balances the support function $x\cdot\nu/(1+s)$, which is nonconstant on every nonspherical member.

\subsection{Main results}
 
Let $n\geq2$, $0<s<1$, and let $E\subset\mathbb R^n$ be a bounded open set with smooth boundary $ \partial E$ and outward unit normal $\nu$. We adopt the following definitions (without factor $s(1-s)$):
\begin{align}
 P_s(E)&=\int_E\int_{\R^n\setminus E}|x-y|^{-n-s}\dd y\dd x,
 \label{ps}\\
 H_s(E,x)&=\pv\int_{\R^n}
 \frac{\chi_{\R^n\setminus E}(y)-\chi_E(y)}{|x-y|^{n+s}}\dd y
 =\frac2s\int_{\partial E}
 \frac{(y-x)\cdot\nu_E(y)}{|x-y|^{n+s}}dH^{n-1}(y).
 \label{eq:boundaryH}
\end{align}
Note that the boundary integral is absolutely convergent for smooth
boundaries. 
Observe that the curvature of every ball is constant, and let us denote from now on $H_s(B_1, \theta)=H_s(\mathbb S^{n-1})$.

For a normal variation with speed $f$, the first variation is
\begin{equation}\label{eq:firstvar}
 \delta P_s(E)[f]=\int_{\partial E}H_s(E,x)f(x)dH^{n-1}(x).
\end{equation}
Thus the associated gradient flow is
\begin{equation}\label{eq:originalflow}
 \partial_tX=-H_s\nu.
\end{equation}

The homogeneity $H_s(\lambda E,\lambda x)=\lambda^{-s}H_s(E,x)$ shows that, at a fixed pole $(x_0,T)$, the rescaled set 
\begin{equation} \widetilde E_q=\ell^{-1}(E_t-x_0), \qquad \text{
with 
}
 q=-\log((T-t)/T),\qquad \ell=(T-t)^{\frac{1}{1+s}}
 \label{eq:fixedrescale}
\end{equation}
has normal velocity
\begin{equation}
\tilde V_\nu=- H_s(\tilde E_q, x)+\frac{1}{1+s} x\cdot\nu_{\tilde E_q}(x).
 \label{eq:defect}
\end{equation}
We fix the notation for the geometric shrinker defect as follows:
\begin{equation}\label{eq:geometricdefect}
 \mathfrak D_s(E,x):=H_s(E,x)-\frac{1}{s+1} x\cdot\nu_E(x),\qquad x\in\partial E.
\end{equation}

  The stationary round shrinker has radius
\begin{equation}
R_s:= \bigl((1+s)H_s(\mathbb S^{n-1})\bigr)^{\frac{1}{1+s}}.
 \label{eq:introRstar}
\end{equation}
 
From now on we fix as reference ball $B_{R_s}$ and we denote it as $B$. 
We introduce the notion of  nearly spherical sets on the ball $B$: 
\begin{definition}[Nearly spherical set]\label{defnearlysph}
For $s\in(0,1)$ let $R_s=((1+s)H_s(\mathbb S^{n-1}))^{1/(1+s)}$.
A set $E_r$ is nearly spherical with respect to $B_{R_s}$ if
\[
 E_r=\{\rho\theta:\theta\in\mathbb S^{n-1},\ 0\leq\rho<r(\theta):=R_s+u(\theta)\},
 \qquad \partial E_r=\{r(\theta)\theta:\theta\in\mathbb S^{n-1}\},
\]
where $0<\eta<1$, $u\in C^{2,\eta}(\mathbb S^{n-1})$ and $\|u\|_{C^{2,\eta}}$ is sufficiently small.
We call $u$ the   radial height. The dimensionless height
$v=u/R_s$ is used only as an auxiliary coordinate when convenient.

We will denote $H_s(E_r,r(\theta)\theta)= H_s(r)(\theta)$.
\end{definition}

We use the smoothness of nonlocal curvature in spherical charts from \cite{CFWlattices}, and the spherical fractional Jacobi spectrum and second-variation formulas from \cite{FFMMM}.

The first nontrivial shape eigenvalue of the rescaled evolution is
\[\mu_2(s)=\frac{n(1-s)-s(1+s)}{(1+s)(n-s)},\]
which vanishes precisely at
\[\sigma=s_*(n)=\frac{\sqrt{n^2+6n+1}-(n+1)}2.\]

We first prove nonlinear stability above the threshold and instability below it. 
For $s>\sigma$, the eigenvalues of the normalized linearization
on harmonics of degrees $k\ge2$ are strictly negative. After removing
translations and dilations, this spectral gap yields nonlinear stability.
The precise statement is Theorem~\ref{thm:dynamicmain}.

\begin{maintheorem}\label{thm:mainA}
For \(s>\sigma\) and \(0<\eta<1\), every sufficiently small
\(C^{2,\eta}\) perturbation of a sphere, up to similarity, evolves
under the unconstrained fractional mean curvature flow to extinction
at a single point in finite time. Its true extinction rescaling
converges exponentially to \(B_{R_s}\) in \(\mathrm h^{2+\eta/2}\)
and in every \(C^k\) norm for normalized times bounded away from zero.
\end{maintheorem}

For $s<\sigma$ we obtain on the other hand that the sphere is dynamically unstable. In particular for $s<\sigma$, we may define the instability index $i(s)$, see \eqref{index},  which is the first harmonic degree with a nonpositive shape eigenvalue. The number of positive eigenvalues counted with multiplicity is $\operatorname{ind}_+(s)$, defined in \eqref{eq:positiveindex}. In particular we get that $i(s)\to +\infty$ as $s\to 0^+$.
At $s=\sigma$, the degree-two eigenvalue vanishes. The strict spectral-gap argument does not decide nonlinear stability at this exponent.

For the precise statement see Theorem \ref{prop:belowthreshold}.

\begin{maintheorem}\label{thm:mainB}
For \(0<s<\sigma\), the round shrinker is nonlinearly Lyapunov unstable modulo translations and dilations in the little-H\"older shape space. 
\end{maintheorem} 

Finally, the third main result of the paper is the existence of convex nonspherical shrinkers for $s$ in a neighborhood of the threshold.  For the precise statement see Theorem \ref{thm:unbalanced} and Theorem \ref{thm:balancedexist}.   
\begin{maintheorem} \label{thm:mainC}
For $n\geq 3$ and each \(1\le p<n/2\) there is a smooth branch of strictly convex, \(O(p)\times O(n-p)\)-invariant homothetic shrinkers bifurcating from the round sphere at \(s=\sigma\). Its nonzero members are nonspherical and occur on both sides of the threshold. 

In even dimensions $n=2p\geq2$, there is also a balanced
branch parametrized by $a$, with
$s(0)=\sigma$, $s'(0)=0,$ $ s''(0)>0.$
Thus its nonzero members lie on the $s>\sigma$ side,
and $s(a)-\sigma$ is quadratic to leading order in $a$.   Moreover for $a\ne0$ sufficiently small these shrinkers are not ellipsoids. 
\end{maintheorem}
 
In dimensions $n\ge3$, we do not classify all nearby shrinkers;
in the plane, a local classification follows from the reduction
(Corollary~\ref{cor:planar}). We conjecture that the nonspherical
branches constructed here are Lyapunov unstable.
An analysis of the reduced evolution on a flow-invariant center
manifold at the threshold may provide a way to address this question.    

\subsection{Notation and functional spaces} \label{sec:notation}
Throughout, $n\geq2$, $0<s<1$, $\mathbb S^{n-1}$ is the unit sphere,
$\nabla=\nabla_{\mathbb S^{n-1}}$ is its tangential gradient, and $I$
denotes the identity operator on the space under consideration.
We write $\langle f\rangle=|\mathbb S^{n-1}|^{-1}\int_{\mathbb S^{n-1}}f dH^{n-1}$ for the spherical average. When integrating on the sphere for sake of brevity we use the notation $\int_{\mathbb S^{n-1}}f(\theta) d\theta$ in place of the more precise one $\int_{\mathbb S^{n-1}}f(\theta) dH^{n-1}(\theta)$.

The symbol $D$ denotes a Fr\'echet derivative, with subscripts specifying
its variable.

For $E\subset\mathbb R^n$, $\nu_E$ denotes the outer unit normal and
$B_R$ the ball of radius $R$ centered at the origin.
We write 
$\sigma=s_*(n)=(\sqrt{n^2+6n+1}-(n+1))/2$.

For $p\in \mathbb{N}$, $p\geq 1$, $O(p)=\{U\in\mathbb R^{p\times p}:U^TU=I_p\}$ is the group of \(p \times p\) orthogonal matrices.
 
The notation $\mathcal L(X,Y)$ means bounded linear operators $X\to Y$.

For $k\in\mathbb N\cup\{0\}$ and $0<\beta<1$, we define the little-H\"older space by
\begin{equation}\label{eq:littleholderdef}
\mathrm h^{k+\beta}= \mathrm h^{k+\beta}(\Sph):=\overline{C^\infty(\Sph)}^{\,C^{k, \beta}(\Sph)}.
\end{equation}
  
We denote by $\mathcal H_k$ the finite-dimensional subspace in $L^2(\mathbb S^{n-1})$ of the spherical harmonics of degree $k$, and with $Y_k^i$ the elements of an orthonormal basis of $\mathcal H_k$. The symbol $P_k$ denotes the orthogonal projection onto $\mathcal H_k$.
In particular, $Q_2=I-P_2$ retains all modes except degree two.
Later $\Pi=P_0+P_1$ projects onto the translation and dilation modes, and $Q=I-\Pi$ removes them;
thus $Q\ne Q_2$. The projections $P_k$ and $\Pi$ have finite rank,
whereas $Q$ and $Q_2$ have finite-dimensional kernels. All these
projections are bounded on the H\"older and little-H\"older spaces
used below.

\subsection{Organization of the paper}
Section~\ref{sec:spectrum} contains the linearization and the spherical
spectrum, and Section~\ref{sec:radial} establishes the radial curvature
estimates. Section~\ref{sec:dynamics} proves dynamical stability above
the threshold and instability below it. We then study stationary
profiles near the critical exponent: Section~\ref{sec:LS} develops
the Lyapunov--Schmidt reduction, and Section~\ref{sec:branches}
constructs the nonspherical shrinking branches.

\subsection*{Acknowledgements} The authors are members of the Gruppo Nazionale per l'Analisi Matematica,
la Probabilit\`a e le loro Applicazioni (GNAMPA)
of the Istituto Nazionale di Alta Matematica (INdAM).
\subsection*{Use of AI tools.} LLMs were used during the preparation of this work to check algebraic identities, and to assist in exploring proof arguments and in revising and reorganizing the manuscript. 
The authors take full responsibility for the mathematical content.

\section{Jacobi operator on the sphere and linearization}
\label{sec:spectrum}
 Throughout the static and dynamical analysis the reference sphere is the
stationary sphere of radius \begin{equation}\label{raggio}
 R_s:= \bigl((1+s)H_s(\mathbb S^{n-1})\bigr)^{1/(1+s)}
\end{equation} and we consider  radial graphs over this sphere, according to Definition \ref{defnearlysph}.

We use the radial height $u$:
\begin{equation} \label{eq:physicalstaticgraph}
 X_{u}(\theta)=(R_s+u(\theta))\theta,\qquad
 E_{r}=\{\rho\theta:0\leq\rho<r(\theta)=R_s+u(\theta)\}.
\end{equation}
 Recall that the stationary shrinker defect is $\mathfrak{D}_s(E,x)= H_s(E, x)-\frac{1}{1+s} x\cdot \nu_E(x)$.

For $\bar y=(R_s+u(\eta))\eta\in\partial E_r$, the outer unit normal is 
$\nu(\bar y)=\frac{(R_s+u(\eta))\eta-\nabla u(\eta)}{\sqrt{ |R_s+u(\eta)|^2+|\nabla u(\eta)|^2}}$ where $\nabla u$ is the tangential gradient of $u$. So, recalling $\nabla u(\eta)\cdot \eta=0$, we get  
\[\bar y\cdot \nu_E(\bar y)=\frac{R_s+u(\eta)}{\sqrt{1+\frac{|\nabla u(\eta)|^2}{|R_s+u(\eta)|^2}}}.\]

Therefore the  stationary shrinker defect  reads as follows: 
\begin{equation}\mathcal F_s[u](\theta):= 
\mathfrak{D}_s(E_r,X_u(\theta))=
 H_s(E_{r},X_{u}(\theta)) -\frac{1}{1+s}\frac{R_s+u(\theta)}{\sqrt{1+\frac{|\nabla u(\theta)|^2}{(R_s+u(\theta))^2}}}.
\label{eq:physicalstaticdefect}
\end{equation}
Thus $\mathcal F_s[u]\equiv0$ is exactly the physical shrinker equation, and $\mathcal F_s[0]\equiv0$
for every $s\in(0,1)$.

For $M=\partial E$ and a scalar normal speed $f:M\to\mathbb R$, the normal derivative of $\Hs$ at the smooth hypersurface $M$ is the Jacobi operator
\begin{align}
 2\pv\int_M\frac{f(x)-f(y)}{|x-y|^{n+s}}d H^{n-1}(y)- 
f(x)\int_M\frac{|\nu(x)-\nu(y)|^2}{|x-y|^{n+s}}dH^{n-1}(y).
 \label{eq:Jacobi}
\end{align}
 
These formulas follow by differentiating the first variation of fractional perimeter; see    the second-variation and evolution formulas in \cite[Section 6]{FFMMM}, \cite[Section 3]{JLMshort}, and \cite{SV} after normalization.

We introduce  the hypersingular Riesz operator on the sphere, with the exact normalization 
\begin{equation}
 \Is u(\theta)=2\pv\int_\Sph
       \frac{u(\theta)-u(\eta)}{|\theta-\eta|^{n+s}}\dd\eta.
 \label{eq:Is}
\end{equation}
Note that \begin{equation}\label{fract}  \int_{\mathbb{S}^{n-1}} u(\theta) \Is u(\theta) dH^{n-1}(\theta) =\int_{\mathbb{S}^{n-1}} \int_{\mathbb{S}^{n-1}} \frac{(u(\theta)-u(\eta))^2}{|\theta-\eta|^{n+s}}  dH^{n-1}(\theta)  dH^{n-1}(\eta):=[u]^2_{\frac{1+s}{2}}   \end{equation}
is the (squared) fractional Gagliardo seminorm of $u$. We will endow  the fractional Sobolev space $H^{\frac{1+s}{2}}(\mathbb{S}^{n-1})$ with the norm given by the square root  of  
\[\|u\|^2_{L^2(\mathbb{S}^{n-1})}+  [u]^2_{\frac{1+s}{2}}. \]

We will denote with $\lambda_k(s)$ the $k$-order eigenvalue of  the operator \eqref{eq:Is}, so there holds that  $ \Is Y_k^i=\lambda_k(s) Y_k^i$, where $Y_k^i$ is an element of an orthonormal basis (in $L^2(\mathbb S^{n-1})$) of  the space $\mathcal{H}_k$ of  spherical harmonics of degree $k$ on the sphere. 
By \cite[Proposition 2.3]{FFMMM} and the references therein,  $\lambda_k(s)>\lambda_{k-1}(s)$ for all $k\geq 1$ and that 
\begin{equation}\label{eig} \lambda_0(s)=0,\quad \lambda_1(s)=s H_s(\mathbb S^{n-1})=s(n-s)\frac{P_s(B_1)}{|\mathbb S^{n-1}|}, \quad \lambda_2(s)= \frac{2n}{n-s}\lambda_1(s). \end{equation}  
Moreover, \begin{equation}\label{conto1} \frac{\lambda_k(s)}{\lambda_1(s)}
 = \sum_{j=0}^{k-1}\frac{((n+s)/2)_j}{((n-s)/2)_j},\qquad k\ge1,
 \end{equation} 
where $(b)_j=b(b+1)\cdots(b+j-1)$ and $(b)_0=1$.

Therefore by \eqref{eq:Jacobi} and the discussion above  the  linearization of the operator  $\mathcal F_s[u]$ defined at \eqref{eq:physicalstaticdefect}  is given by 
\begin{equation}
 D_u\mathcal F_s[0]=\frac{1}{R_s^{1+s}}\Is-\frac{sH_s(B)}{R_s}I-\frac{1}{1+s}I= \frac{s}{(s+1)\lambda_1(s)}\Is-I:= \mathcal J_s
 \label{eq:As}
\end{equation}
where we used the fact that $R_s^{1+s}= \frac{1+s}{s}\lambda_1(s)$ and $H_s(B)=
H_s(\mathbb S^{n-1})R_s^{-s} $. 
Thus, in fixed-pole rescaled time $q$, the linearized evolution is \[\partial_q u=-\mathcal J_su.\] On harmonics of degree $k$,
\begin{equation}
-\mathcal J_s Y_k=\mu_k(s)Y_k,\qquad
 \mu_k(s)=1-\frac{s}{1+s}\frac{\lambda_k(s)}{\lambda_1(s)}.
 \label{eq:muk}
\end{equation}
The degree-zero and degree-one eigenvalues are $1$ and $\frac{1}{1+s}$, respectively. They represent changes of scale/extinction parameter and translations in fixed-pole rescaled coordinates.

Note that by \eqref{eig} we have an explicit formula  for $\mu_2(s)$:
\begin{equation}\label{mu2}\mu_2(s)=1- \frac{2ns}{(n-s)(1+s)}=\frac{n-(n+1)s-s^2}{(n-s)(s+1)}.\end{equation}
 
As $s\uparrow1$, $\mu_2(s)\to-\frac{1}{n-1}$ and there is a threshold
\[\sigma:=s_*(n)= \frac{\sqrt{n^2+6n+1}-n-1}{2}\] such that for $s>\sigma $, it holds that $\mu_2(s)<0$, and so  $\mu_k(s)<0$ for $k>2$,  whereas $\mu_2(\sigma)=0$ and $\mu_2(s)>0$ for $s\in (0,\sigma)$. 

For $s<\sigma$ define the first nonunstable harmonic degree  
\begin{equation}\label{index} i(s)=\min \{k>2, \mu_k(s)\leq 0\}.\end{equation} 
The shape instability index counted with multiplicity is
\begin{equation}\label{eq:positiveindex}
 \operatorname{ind}_+(s)=\sum_{k=2}^{i(s)-1}\dim\mathcal H_k.
\end{equation}
Here $\dim\mathcal H_k=\binom{n+k-1}{k}-\binom{n+k-3}{k-2}$
for $k\ge2$. We get that $i(s)$  is finite for each $s>0$, and $i(s)\to +\infty$ as $s\to 0^+$,  since by \eqref{conto1} we get
\begin{equation}
 \mu_k(s)=1-\frac{s}{s+1}\frac{\lambda_k(s)}{\lambda_1(s)}
 =1-\frac{s}{s+1}\sum_{j=0}^{k-1}\frac{((n+s)/2)_j}{((n-s)/2)_j},\qquad k\ge1.
 \label{eq:pochhammer}
\end{equation}

  Let $\Delta_j$, $j\ge0$, be a smooth dyadic partition
in spherical harmonic degree, with $j=0$ containing the low degrees.
For every positive noninteger $b$ one has
\begin{equation}\label{eq:dyadicHolder}
 \|f\|_{C^b}\asymp\sup_{j\ge0}2^{jb}\|\Delta_jf\|_\infty;
 \qquad
 f\in\mathrm h^b\ \Longleftrightarrow\
 2^{jb}\|\Delta_jf\|_\infty\longrightarrow0.
\end{equation}
These are the usual H\"older--Zygmund characterization and its little
subspace at noninteger exponents. The smooth spherical multiplier
localization theorem gives the following concrete input: a multiplier
supported on degrees comparable to $2^j$, whose rescaled symbol has
bounded $C^N$ norm, has kernel bounded by
\[
 C_N2^{j(n-1)}(1+2^j\operatorname{dist}(\theta,\eta))^{-N},
\]
for $N$ large enough, with a constant controlled by finitely many
symbol derivatives. Its kernel therefore has bounded $L^1$ norm.
This localization statement follows from \cite[Theorem 3.5]{NPW},
using a smooth annular symbol; on the circle it is the elementary
Fourier multiplier estimate obtained by summation by parts.
Formula \eqref{eq:dyadicHolder} follows from the corresponding
Bernstein inequalities and smooth approximation estimates. In
particular a symbol $a(t)$ satisfying
$|\partial_t^N a(t)|\le C_N(1+t)^{d-N}$ induces a bounded map
$\mathrm h^{b+d}\to\mathrm h^b$ whenever both exponents are positive
and noninteger: apply the localized kernel estimate on each block,
then \eqref{eq:dyadicHolder}. All constants are uniform on bounded
sets of these symbol seminorms.

\begin{lemma} \label{lem:spectralresolvent}
Fix a compact interval $I\Subset(0,1)$ sufficiently close to $\sigma$.
Set $\mathscr L_s=s\mathcal I_s/\lambda_1(s)-(1+s)I$, the normalized static linearization. The operator $Q_2\mathscr L_sQ_2$ is an isomorphism from
$Q_2\mathrm h^{b+1+s}$ onto $Q_2\mathrm h^b$, uniformly for $s\in I$
when the exponents stay in compact noninteger ranges. For every fixed $s\in(0,1)$,
$-\mathcal I_s$ generates a strongly continuous analytic semigroup
on $\mathrm h^b$, with domain $\mathrm h^{b+1+s}$.
\end{lemma}
\begin{proof}
The ratio formula
\eqref{conto1} determines all eigenvalues. The smooth extension of
the eigenvalue sequence at large $t$ is
\begin{equation}\label{eq:contoGamma}
\begin{gathered}
 \lambda_t(s)=C_{n,s}\left[
 \frac{\Gamma(t+(n+s)/2)}{\Gamma(t+(n-2-s)/2)}
 -\frac{\Gamma((n+s)/2)}{\Gamma((n-2-s)/2)}\right],\\
 C_{n,s}=\frac{2^{1-s}\pi^{(n-1)/2}\Gamma((1-s)/2)}
 {(1+s)\Gamma((n+s)/2)}>0.
\end{gathered}
\end{equation}
It agrees with $\lambda_k(s)$ at $t$ integers, by taking consecutive
differences and using $\lambda_0=0$, $\lambda_1=sH_s(\mathbb S^{n-1})$;
the value of the latter also follows by a beta integral in
\eqref{eq:boundaryH}. Stirling's formula with differentiated remainder
gives, locally uniformly in $s$,
\[
 \lambda_t(s)\asymp t^{1+s},\qquad
 |\partial_t^N\lambda_t(s)|\le C_N t^{1+s-N}\quad(t\ge t_0).
\]
Thus the inverse complementary static symbol
$[s\lambda_t(s)/\lambda_1(s)-(1+s)]^{-1}$ has order $-1-s$
at high degrees. At the remaining finitely many degrees, all
eigenvalues except degree two are uniformly separated from zero
after shrinking $I$. Setting the degree-two multiplier to zero and
applying the block estimate above proves the inverse bounds and
surjectivity by the convergent dyadic expansion in the little spaces.

The resolvent argument below applies to every fixed $s\in(0,1)$,
independently of the restriction on $I$. Choose
$\varepsilon\in(0,\pi/2)$ and
$|\arg z|\le\pi-\varepsilon$. Positivity of $\lambda_k(s)$ gives
$|z+\lambda_k(s)|\ge c_\varepsilon(|z|+\lambda_k(s))$.
On the $j$th high-frequency block the rescaled derivatives of
$(z+\lambda_t(s))^{-1}$ are bounded by
$C_N(|z|+2^{j(1+s)})^{-1}$.
The block estimate consequently yields, for $|z|\ge1$,
\begin{equation}\label{eq:resolventestimate}
 |z|\|f\|_{\mathrm h^b}+\|f\|_{\mathrm h^{b+1+s}}
 \le C_\varepsilon\|(z+\mathcal I_s)f\|_{\mathrm h^b}.
\end{equation}
The same multipliers construct an onto inverse. They preserve the
vanishing dyadic tails in \eqref{eq:dyadicHolder}, so this is a
resolvent on the little spaces, not just on the large H\"older spaces.
Closedness and the domain characterization follow from this estimate.
The analytic-generator criterion and density of smooth functions give
the stated strongly continuous analytic semigroup. These estimates
also imply compactness of the resolvent and elliptic bootstrapping.
\end{proof}

\section{Fractional perimeter and curvature in radial coordinates}
\label{sec:radial}
We recall exact formulas for fractional perimeter and curvature
in radial coordinates. For fixed $s$, smoothness of the curvature
operator in these coordinates is established in
\cite[Theorem~2.1]{CFWlattices}. Here we give the higher-regularity
estimates needed below and compute the second radial derivative
of curvature at the sphere.

For a positive $r\in C^2(\mathbb S^{n-1})$, let
$E_r=\{\rho\theta:0\le\rho<r(\theta)\}$.
The radial perimeter representation follows from
\cite[proof of Theorem~2.1, Step~1, between (2.19) and (2.20)]{FFMMM};
see also \cite{CNball}. By homogeneity it applies to any positive
reference radius. With our normalization it reads 
\begin{equation}\label{eq:radialexactnewper} 
 P_s(E_r)
 =\frac{H_s(\mathbb S^{n-1})}{n-s}
 \int_{\mathbb S^{n-1} }r(\theta)^{n-s}\,\mathrm d\theta
+\frac12
 \int_{\mathbb S^{n-1}}\int_{ \mathbb S^{n-1}} 
 \int_{r(\eta)}^{r(\theta)}
 \int_{r(\eta)}^{r(\theta)}
 \frac{t^{n-1}\rho^{n-1}}
 {|t\theta-\rho\eta|^{n+s}}\,\mathrm d\rho\,\mathrm dt\,\mathrm d\eta\,\mathrm d\theta.
\end{equation}
Both inner integrals are oriented from $r(\eta)$ to $r(\theta)$.
Their signs therefore cancel, and the second term is nonnegative
and finite. The curvature is given by 
\begin{equation}\label{eq:radialexactnew}
  H_s(r)(\theta)
=H_s(E_r,r(\theta)\theta)=\frac2s\int_{\mathbb S^{n-1}}\frac{r(\eta)^{n-2}\,T_r(\theta,\eta)}{W_r(\theta,\eta)^{(n+s)/2}|\theta-\eta|^{n+s}}\,d\eta,
\end{equation}
where 
\begin{align}
 W_r(\theta,\eta)&=\frac{(r(\theta)-r(\eta))^2}{|\theta-\eta|^2}+r(\theta)r(\eta)  ,\label{eq:radialQnew}\\
 \label{eq:radialRnew}
 \mathcal R r(\theta,\eta)&=r(\eta)-r(\theta)-(\eta-\theta)\cdot\nabla r(\eta)\\T_r(\theta,\eta)&=r(\theta)\mathcal R r(\theta,\eta)+(r(\theta)-r(\eta))^2+\tfrac12 r(\theta)r(\eta) |\theta-\eta|^2.\label{eq:radialTnew}
\end{align}

\begin{lemma} \label{lem:diagonalnew}
Fix $0<s_-<s_+<1$, an integer $k\ge0$, and $M,c>0$. Suppose $r\in C^{k+2}(\mathbb S^{n-1})$ satisfies $r\ge c$ and $\|r\|_{C^{k+2}}\le M$. For $|\theta-\eta|\le1$ one has
\begin{equation}\label{eq:diagonalnew}
 W_r(\theta, \eta)\ge c^2,\qquad |\nabla_{\rm rot}^{\le k}T_r(\theta, \eta)|\le C |\theta-\eta|^2,
\end{equation}
where angular derivatives mean simultaneous derivatives under rotations of both variables. For directions $h_1,\ldots,h_m\in C^{k+2}$, $m\ge1$, the same estimates hold for $D_r^mT_r[h_1,\ldots,h_m]$, with right-hand side $C|\theta-\eta|^2\prod_i\|h_i\|_{C^{k+2}}$. Furthermore, derivatives of $ W_r(\theta, \eta) $ of any fixed finite order in $r$ and up to $k$ simultaneous angular derivatives are bounded, with the corresponding products of direction norms.
\end{lemma}
\begin{proof}
  Taylor's theorem   together with $\nabla r(\eta)\cdot\eta=0$, gives
$|\mathcal R r(\theta,\eta)|\le C\|r\|_{C^2}|\eta-\theta|^2$.
The same estimate holds after applying a simultaneous rotation derivative: such derivatives commute with the geometric remainder and replace $r$ by its infinitesimal rotation derivatives. Since $|r(\theta)-r(\eta)|\le C|\eta-\theta|$, every term of $T_r=r(\theta)\mathcal R r+(r(\theta)-r(\eta))^2+ r(\theta)r(\eta) |\eta-\theta|^2/2$ is $O(|\eta-\theta|^2)$. The expression is polynomial in $r$ and its first derivatives, so differentiation with respect to $r$ preserves this property. Finally
$W_r=r(\theta)r(\eta)+((r(\theta)-r(\eta))/|\eta-\theta|)^2\ge c^2$,
and the difference quotient $(r(\theta)-r(\eta))/|\theta-\eta|$ and its variations have the stated bounds. These arguments also apply to simultaneous rotation derivatives, because $|\theta-\eta|$ is rotation invariant.
\end{proof}

\begin{proposition} \label{prop:cutoffnew}
Let $s\in[s_-,s_+]\Subset(0,1)$ and let $r$ range over a bounded positive subset of $C^{k+2}(\mathbb S^{n-1})$. Define $ H_{s,\varepsilon}(r)(\theta)$ by restricting \eqref{eq:radialexactnew} to $\eta$ such that $|\theta-\eta|>\varepsilon$. For every fixed nonnegative integer $m$, the $m$-th Fr\'echet derivative in $r$ of the truncated integral converges, as $\varepsilon\downarrow 0$, uniformly in $\theta$ and in $C^k$ under simultaneous rotation differentiation. More precisely,
\begin{equation}\label{eq:cutoffrate}
 \big\|D_r^m H_{s,\varepsilon}(r)[h_1,\ldots,h_m]-D_r^mH_s(r)[h_1,\ldots,h_m]\big\|_{C^k}
 \le C\varepsilon^{1-s_+}\prod_{i=1}^m\|h_i\|_{C^{k+2}}.
\end{equation}
 Consequently $r\mapsto H_s(r)$ is $C^\infty$ from a positive open subset of $C^{k+2}$ into $C^k$, and differentiation under the untruncated integral is legitimate. Moreover, for every fixed $j\ge0$, the mixed derivatives $\partial_s^jD_r^m H_s$ exist and satisfy the cutoff estimate
\begin{equation}\label{eq:cutofflogs}
 C_{j,m,k}\varepsilon^{1-s_+}(1+|\log\varepsilon|^j)\prod_i\|h_i\|_{C^{k+2}}.
\end{equation}
Thus $(s,r)\mapsto H_s(r)$ is jointly $C^\infty$ with values in $C^k$ on these positive $C^{k+2}$ sets.
\end{proposition}
\begin{proof}
  By Lemma~\ref{lem:diagonalnew}, $W_r\ge c^2$, and all the indicated derivatives of $W_r$ are uniformly bounded. The integrand equals
\[
 \frac2s\,|\theta-\eta|^{-n-s}\,r(\eta)^{n-2}W_r^{-(n+s)/2}T_r.
\]
Repeated product and chain rules, including arbitrary $r$-variations, therefore give the \emph{combined} pointwise estimate
\begin{equation}\label{eq:integrandboundnew}
 \left|\nabla_{\rm rot}^{\le k}D_r^m\left[\frac2s\,|\theta-\eta|^{-n-s}\,r(\eta)^{n-2}W_r^{-(n+s)/2}T_r\right][h_1,\ldots,h_m]\right|
 \le C|\theta-\eta|^{2-n-s}\prod_i\|h_i\|_{C^{k+2}}.
\end{equation}
 
Near the diagonal, spherical polar coordinates give a surface
measure bounded by $C_n\rho^{n-2}\,\mathrm d\rho$, with
$\rho=|\theta-\eta|<1$. Hence the omitted diagonal cap is bounded by
\[
 C\int_0^\varepsilon\rho^{2-n-s}\rho^{n-2}\,d\rho
 =\frac{C}{1-s}\varepsilon^{1-s}\le C'\varepsilon^{1-s_+}.
\]
Simultaneous rotations preserve $|\theta-\eta|$ and the cutoff domain; consequently rotation derivatives pass under the truncated integral without boundary terms. Equivalence of rotation and coordinate derivatives on the compact sphere yields the $C^k$ estimate. Taylor's theorem for the coefficient functions and the same integrable majorant justify Fr\'echet differentiability, including continuity of every derivative. Finally, differentiating $|\theta-\eta|^{-n-s}$ in $s$ introduces $(-\log |\theta-\eta|)^j$, while derivatives of $s^{-1}$ and $W_r^{-(n+s)/2}$ introduce only bounded coefficient factors and powers of $\log W_r$. Integrating $\rho^{-s}(1+|\log\rho|^j)$ gives \eqref{eq:cutofflogs}, and the same dominated-convergence argument proves the claimed joint smoothness.
\end{proof}

Using the notation introduced in \eqref{eq:radialQnew} and \eqref{eq:radialRnew}, we define a linear operator in $v$, with coefficients determined by $r$,   as follows:
\begin{align}
 \mathcal A_s(r)[v](\theta)
 =\frac2s\int_\Sph\frac{r(\eta)^{n-2}}{W_r^{\frac{n+s}{2}}|\theta-\eta|^{n+s}}
 \bigg\{&r(\theta)\mathcal R v(\theta,\eta)
       +\bigl(r(\eta)-r(\theta)\bigr)\bigl(v(\eta)-v(\theta)\bigr)\notag\\
       &+\frac{r(\theta)v(\eta)}2|\theta-\eta|^2\bigg\}\dd\eta.
 \label{eq:quasilinearradial}
\end{align}
This integral is absolutely convergent whenever $r$ is $C^1$ and positive and $v\in C^{1+\gamma}$ with $\gamma>s$. Indeed recall that   for  $r$ bounded below by a positive constant and bounded in $C^1$, $W_r$ is bounded above and below uniformly. 
Moreover  for $v\in C^{1+\gamma}$, $0<\gamma<1$,   $\mathcal R v(\theta,\eta) $ is  $O(|\theta-\eta|^{1+\gamma})$; for $v\in C^2$ it is $O(|\theta-\eta|^2)$.

\begin{lemma} 
\label{lem:exactradialsplit}
For positive $C^2$ radial functions $r(\theta)=R_s+u(\theta)$,
\begin{equation}
H_s(r)=\mathcal A_s(r)[r]
     =\mathcal A_s(r)[u]+\mathcal B_s(r),
 \label{eq:Hsplit}
\end{equation}
where
\begin{equation}
 \mathcal B_s(r)=\mathcal A_s(r)[R_s]
 =\frac{R_s}{s}r(\theta)\int_\Sph
          \frac{r(\eta)^{n-2} }{{W_r^{\frac{n+s}{2}}|\theta-\eta|^{n+s-2}}}\dd\eta.
 \label{eq:radiallowerterm}
\end{equation}
\end{lemma}
\begin{proof}
For $y=r(\eta)\eta$ one has
\[
 \nu(y)dH^{n-1}(y)
   =\bigl(r(\eta)^{n-1}\eta-r(\eta)^{n-2}\nabla r(\eta)\bigr)\dd\eta.
\]
Insert this into the boundary formula \eqref{eq:boundaryH}. Since
$\theta\cdot\eta=1-|\theta-\eta|^2/2$ and $\eta\cdot\nabla r(\eta)=0$, its numerator becomes
\begin{align*}
r(\eta)^{n-2}\bigl(r(\theta)[r(\eta)-r(\theta)+(\theta-\eta)\cdot\nabla r(\eta)]+(r(\eta)-r(\theta))^2+\tfrac12r(\theta)r(\eta)|\theta-\eta|^2\bigr).
\end{align*}
This proves $H_s(r)=\mathcal A_s(r)[r]$; linearity in the bracketed argument gives the remainder of \eqref{eq:Hsplit}. Notice that $\mathcal B_s$ contains only values of its argument and first difference quotients in its coefficients.
\end{proof}

Fix   $\delta\in(0,1)$ and choose
\begin{equation}
 \max\{s,\delta\}<\gamma<\min\{1,s+\delta\},
 \qquad \gamma-s<\delta_0<\delta.
 \label{eq:Holderchoices}
\end{equation}

\begin{lemma}
\label{lem:radialmapping}
Let $j\geq1$ be an integer. In a bounded positive radial neighborhood in $\mathrm h^{j+1+\delta_0}$,
\begin{align}
 \mathcal A_s&\in C^\infty\!\left(
     \mathrm h^{j+1+\delta_0},
     \mathcal L(\mathrm h^{j+1+\gamma},\mathrm h^{j+\gamma-s})\right),
 \label{eq:Amapping}\\
 \mathcal B_s&\in C^\infty\!\left(\mathrm h^{j+1+\delta_0}, \mathrm h^{j+\delta_0}\right).
 \label{eq:Bmapping}
\end{align}
More explicitly, for every $k\geq0$ there is $C_k$ such that
\begin{align}
 \|D^k\mathcal A_s(r)[w_1,\ldots,w_k]v\|_{C^{j+\gamma-s}}
 &\leq C_k\|v\|_{C^{j+1+\gamma}}
               \prod_{i=1}^k\|w_i\|_{C^{j+1+\delta_0}},
\label{eq:Adifferentialbound}\\
 \|D^k\mathcal B_s(r)[w_1,\ldots,w_k]\|_{C^{j+\delta_0}}
 &\leq C_k\prod_{i=1}^k\|w_i\|_{C^{j+1+\delta_0}}.
 \label{eq:Bdifferentialbound}
\end{align}
Constants depend on a positive lower bound for $r$ and an upper bound for its indicated norm. The statements include local Lipschitz estimates for these derivatives.
\end{lemma}
\begin{proof}
We give the cancellation and H\"older estimates, including the potentially singular terms. Rotation vector fields span the tangent bundle of the sphere, and their iterates give equivalent $C^{j+\vartheta}$ norms. To differentiate the value at $\theta$ in a rotation direction, rotate both $\theta$ and the integration variable $\eta$. The surface measure $\mathcal H^{n-1}$ and  $|\theta-\eta|$ are invariant. Differentiation therefore falls only on radial functions, their gradients, and their difference quotients.

After at most $j$ such differentiations, every coefficient in \eqref{eq:quasilinearradial} is a smooth function of uniformly bounded quantities of the form
\[
 r(\theta),\ r(\eta),\quad
 \frac{r(\eta)-r(\theta)}{|\theta-\eta|},
 \quad\text{and their simultaneous rotation derivatives.}
\]
The same assertion holds after parameter derivatives, with $r$ replaced in some factors by $w_i$. In particular
\[
 \frac{|w(\eta)-w(\theta)|}{|\theta-\eta|}\leq C\|w\|_{C^1}.
\]
If a simultaneous rotation has size $h\leq1$, a coefficient with up to $j$ derivatives changes by at most
$Ch^{\delta_0}$ times the relevant product of $C^{j+1+\delta_0}$ norms. This follows from the $C^{\delta_0}$ modulus of the first derivatives of the $j$th rotation derivatives. The positive lower bound for $ W_r$ makes the same estimates valid for every differentiated denominator.

It remains to estimate the three numerator types. For the principal type, the radial absolute-integrability estimate is, letting $d=|\theta-\eta|$, 
\begin{equation}
 \int_0^1 d^{1+\gamma}d^{-1-(1+s)}\dd d
       =\int_0^1 d^{\gamma-s-1}\dd d<\infty.
 \label{eq:principalintegrability}
\end{equation}
To estimate a difference of output values at angular distance $h$, first change variables by a rotation. Split the integral into $d\leq h$ and $d>h$. The principal Taylor remainder contributes, on the first region, at most $C h^{\gamma-s}$. On the second region use
\[
 \|v\circ U-v\|_{C^1}\leq C h^\gamma\|v\|_{C^{1+\gamma}}
\]
for a rotation $U$ of size $h$; hence its Taylor remainder is bounded by
$C h^\gamma |\theta-\eta|\|v\|_{C^{1+\gamma}}$. Integration yields
\begin{equation}
 h^\gamma\int_h^1 d^{-1-s}\dd d
        \leq C h^{\gamma-s}.
 \label{eq:farHolderprincipal}
\end{equation}
The change of coefficient contributes $O(h^{\delta_0})$ times the finite integral \eqref{eq:principalintegrability}, which is bounded by $O(h^{\gamma-s})$ because $\delta_0>\gamma-s$.

For the other two numerator types, factor out $|\theta-\eta|^2$. The remaining coefficients contain either two first difference quotients or a value of $v$. Their integrability is
\begin{equation}
 \int_0^1 d^2d^{-1-(1+s)}\dd d
      =\int_0^1 d^{-s}\dd d<\infty.
 \label{eq:lowerintegrability}
\end{equation}
They have simultaneous rotation modulus $O(h^{\delta_0})$ (or the stronger modulus supplied by $v$). This proves the $C^{\gamma-s}$ estimate. Applying the same argument after $j$ rotation derivatives proves \eqref{eq:Adifferentialbound}; terms with fewer derivatives are bounded by the stronger available norms. For   $\mathcal B_s$ only the last type occurs, with  $|\theta-\eta|^2$,  so the full exponent $\delta_0$ is retained, proving \eqref{eq:Bdifferentialbound}.

Finally, Taylor expansion of the smooth coefficient functions on their uniformly positive domain gives the stated derivatives and a quadratic Taylor remainder in the coefficient norm. This proves Fr\'echet smoothness. Smooth approximation of $r,v,w_i$ and the same estimates show that the maps and derivatives take values in the little H\"older spaces; smooth positive radial data give smooth outputs. This proves \eqref{eq:Amapping}--\eqref{eq:Bmapping}.
\end{proof}

Finally we compute the second derivative of curvature in radial coordinates at radius $r=R_s$ (so at $u=0$).
\begin{proposition} 
\label{prop:second-radial-curvature}
Let $n\ge2$ and $s\in(0,1)$. 
For $u,v\in C^\infty(\mathbb S^{n-1})$, the second radial
derivative at the constant function $r\equiv R_s$ is
\begin{align}
 D_r^2 H_s(R_s)[u,v]
 =R_s^{-s-2}\Bigl\{
 &-\frac{n+s}{2}\bigl(u\mathcal I_sv+v\mathcal I_su\bigr)
 +\frac{n-2-s}{2}\mathcal I_s(uv)\notag\\
 &+s(1+s)H_s(\mathbb S^{n-1})uv
 \Bigr\}.
 \label{eq:second-radial-curvature-Rs}
\end{align}
Here the derivative  is the mixed derivative of
$H_s(R_s+\varepsilon u+\delta v)$ at
$(\varepsilon,\delta)=(0,0)$.
\end{proposition}

\begin{proof}
We first work at radius one and then use homogeneity.
All integrals without an indicated domain are over $\mathbb S^{n-1}$.

For $h\in C^\infty(\mathbb S^{n-1})$ and $\varepsilon$ sufficiently small,
we claim that
\begin{align}
 P_s(E_{1+\varepsilon h})
 ={}&  P_s(E_1)+\varepsilon H_s(\mathbb S^{n-1})\int h d\theta 
 \label{eq:perimeter-cubic-expansion}\\
 &+\frac{\varepsilon^2}{2}\int h\mathcal I_s h\,\mathrm d\theta
 +\frac{\varepsilon^2}{2}(n-1-s)H_s(\mathbb S^{n-1})
       \int h^2\,\mathrm d\theta\notag\\
 &+\frac{\varepsilon^3}{4}(n-2-s)
       \int h^2\mathcal I_s h\,\mathrm d\theta\notag\\
 &+\frac{\varepsilon^3}{6}(n-1-s)(n-2-s)
       H_s(\mathbb S^{n-1})\int h^3\,\mathrm d\theta
       +O(\varepsilon^4).
\notag
\end{align}

This expansion follows from the exact perimeter formula \eqref{eq:radialexactnewper}.
The local term is expanded by the ordinary Taylor formula:
\begin{align*}
 \frac{H_s(\mathbb S^{n-1})}{n-s}\int(1+\varepsilon h)^{n-s}
 ={}&P_s(E_1)+\varepsilon H_s(\mathbb S^{n-1})\int h
 +\frac{\varepsilon^2}{2}(n-1-s)H_s(\mathbb S^{n-1})\int h^2\\
 &+\frac{\varepsilon^3}{6}(n-1-s)(n-2-s)H_s(\mathbb S^{n-1})\int h^3
 +O(\varepsilon^4).
\end{align*}
For the nonlocal term, set $U=h(\theta)$, $V=h(\eta)$. Parametrize the two
oriented intervals by
\[
 t=1+\varepsilon\bigl(V+a(U-V)\bigr),\qquad
 \rho=1+\varepsilon\bigl(V+b(U-V)\bigr),
 \qquad a,b\in[0,1].
\]
The product of their Jacobians is $\varepsilon^2(U-V)^2$.
Moreover,
\[
 |t\theta-\rho\eta|^2
 =|\theta-\eta|^2
   \left(t\rho+\frac{(t-\rho)^2}{|\theta-\eta|^2}\right).
\]
Since $|U-V|\leq\|h\|_{C^1}|\theta- \eta|$, the quantity in parentheses
is bounded away from zero and has bounded derivatives in
$\varepsilon$, uniformly in $\theta\ne\eta$ and $a,b\in[0,1]$.
Thus Taylor expansion under the $a,b$ integrals gives
\begin{align*}
 &\int_{1+\varepsilon V}^{1+\varepsilon U}
 \int_{1+\varepsilon V}^{1+\varepsilon U}
 \frac{t^{n-1}\rho^{n-1}}{|t\theta-\rho\eta|^{n+s}}
 \,\mathrm d\rho\,\mathrm dt\\
 &\qquad=\frac{\varepsilon^2(U-V)^2}{|\theta- \eta|^{n+s}}
+\frac{n-2-s
}{2}\frac{\varepsilon^3(U+V)(U-V)^2}{|\theta- \eta|^{n+s}}
 +O\bigl(\varepsilon^4|\theta- \eta|^{2-n-s}\bigr).
\end{align*}
Indeed, the derivative of the rescaled kernel at $\varepsilon=0$
is $\frac{n-2-s
}{2}\bigl(2V+(a+b)(U-V)\bigr)$, whose average over
$[0,1]^2$ is $\frac{n-2-s}{2}(U+V)$. The remainder is integrable over
$\mathbb S^{n-1}\times\mathbb S^{n-1}$, since its radial bound near the
 diagonal is a constant times $\int_0^1 \rho^{-s}\,\mathrm d\rho<\infty$.
Finally, the normalization and symmetry of $\mathcal I_s$ give
\begin{align*}
 \int h\mathcal I_s h
 &=\int\int
 \frac{(U-V)^2}{|\theta-\eta|^{n+s}}\,\mathrm d\eta\,\mathrm d\theta,\\
 \int h^2\mathcal I_s h
 &=\int\int
 \frac{(U+V)(U-V)^2}{|\theta-\eta|^{n+s}}\,\mathrm d\eta\,\mathrm d\theta.
\end{align*}
This proves 
\eqref{eq:perimeter-cubic-expansion}.
The estimates are uniform when $h$ ranges over a bounded subset
of any fixed finite-dimensional space of smooth functions;
parameter differentiation is justified by the same integrable
majorant.

The coefficients of $\varepsilon^2$ and $\varepsilon^3$ in
\eqref{eq:perimeter-cubic-expansion} are, respectively,
$D^2 P_s(E_1)[h,h]/2$ and
$D^3  P_s(E_1)[h,h,h]/6$.
Therefore it  yields
\begin{equation}
 D^2 P_s(E_1)[u,z]
 =\int z\mathcal I_su+(n-1-s)H_s(\mathbb S^{n-1})\int uz
 \label{eq:second-perimeter-radial}
\end{equation}
and 
\begin{align}
 D^3\  P_s(E_1)[u,v,z]
 ={}&\frac{n-2-s}{2}\int
 (uv\mathcal I_sz+uz\mathcal I_sv+vz\mathcal I_su)
 \notag\\
 &+(n-1-s)(n-2-s)H_s(\mathbb S^{n-1})\int uvz.
 \label{eq:third-perimeter-radial}
\end{align}

The radial first-variation formula is
\[
 D P_s(E_r)[z]:=\left.\frac{d}{d\varepsilon}\right|_{\varepsilon=0}
P_s(E_{r+\varepsilon z})=\int_{\mathbb S^{n-1}}
H_s(r)(\theta)\,z(\theta)\,r(\theta)^{n-1}\,d\theta.
\] 
 
Indeed, 
the normal component of the radial velocity $V(r(\theta)\theta)=z(\theta)\theta$ multiplied
by the surface Jacobian $r^{n-2}\sqrt{r^2+|\nabla_{\mathbb S}r|^2}
$ equals $zr^{n-1}$.
Differentiating at $r=1$ gives
\[
 D^2 P_s(E_1)[u,z]
 =\int    D_rH_s(1)[u] z +(n-1)H_s(\mathbb S^{n-1})u z.
\]
Comparing  with \eqref{eq:second-perimeter-radial}, for arbitrary
$z$,  we obtain 
\[
 D_rH_s(1)[u]=\mathcal I_su-sH_s(\mathbb S^{n-1})u.
\]
Differentiating the first-variation formula twice gives
\begin{align*}
 D^3  P_s(E_1)[u,v,z]
 =\int z\Bigl\{
 &D_r^2H_s(1)[u,v]
 +(n-1)\bigl(uD_rH_s(1)[v]+vD_rH_s(1)[u]\bigr)\\
 &+(n-1)(n-2)H_s(\mathbb S^{n-1})uv\Bigr\}.
\end{align*}
By self-adjointness, $\int uv\mathcal I_sz=\int z\mathcal I_s(uv)$.
We can therefore compare this identity with
\eqref{eq:third-perimeter-radial} and use the arbitrariness of $z$.
Substituting the first derivative gives
\begin{align*}
 D_r^2H_s(1)[u,v]
 ={}&\left(\frac{n-2-s}{2}-(n-1)\right)
       (u\mathcal I_sv+v\mathcal I_su)
 +\frac{n-2-s}{2}\mathcal I_s(uv)\\
 &+\Bigl((n-1-s)(n-2-s)+2s(n-1)
              -(n-1)(n-2)\Bigr)H_s(\mathbb S^{n-1})uv.
\end{align*}
Observe that 
\[
 \frac{n-2-s}{2}-(n-1)=-\frac{n+s}{2},\qquad
 (n-1-s)(n-2-s)+2s(n-1)-(n-1)(n-2)=s(1+s).
\]
This proves the asserted formula at radius one.

Finally the homogeneity of fractional mean curvature gives
\[
 H_s(R_s+\varepsilon u+\delta v)
 =R_s^{-s}H_s\left(1+\frac{\varepsilon}{R_s}u
                         +\frac{\delta}{R_s}v\right).
\]
Taking the mixed derivative in $\varepsilon$ and $\delta$
produces the factor $R_s^{-s-2}$ and proves
\eqref{eq:second-radial-curvature-Rs}.
\end{proof}

\section{Dynamical stability and instability  of the sphere}
\label{sec:dynamics}

Throughout this section, $s\in(0,1)$ is fixed. We use the radial height $u$ and
the stationary defect $\mathcal F_s[u]$ from
\eqref{eq:physicalstaticdefect}, with $r=R_s+u$. The curvature
$H_s(r)$ and its splitting
$H_s(r)=\mathcal A_s(r)[u]+\mathcal B_s(r)$ are those of
Section~\ref{sec:radial}.

We first establish the relation between the shrinker equation and
the rescaled evolution. We then choose a time-dependent translation
and dilation so that the radial perturbation has vanishing
degree-zero and degree-one components. This yields an exact
evolution equation for the normalized profile, to which we apply
quasilinear parabolic theory. We subsequently reconstruct the
original fractional mean curvature flow and prove convergence
under rescaling about its extinction point and extinction time.
The same formulation also yields instability of the sphere
for $s<\sigma=s_*(n)$.

For the radial parametrization $X_u(\theta)=(R_s+u(\theta))\theta$,
one has $\theta\cdot\nu=1/\sqrt{1+\frac{|\nabla u|^2}{(R_s+u)^2}}$ and
$X_u\cdot\nu=(R_s+u)/\sqrt{1+\frac{|\nabla u|^2}{(R_s+u)^2}}$. 
In the fixed-pole rescaled time $q=-\log((T-t)/T)$, the normal
velocity is $u_q/\sqrt{1+\frac{|\nabla u|^2}{(R_s+u)^2}}$. Since this velocity equals
$-\mathcal F_s[u]$, by \eqref{eq:defect}, the exact radial
evolution is
\begin{equation}
 u_q=-\sqrt{1+\frac{|\nabla u|^2}{(R_s+u)^2}}\mathcal F_s[u]
     =-\sqrt{1+\frac{|\nabla u|^2}{(R_s+u)^2}} H_s(R_s+u)+\frac{1}{1+s}(R_s+u).
\label{eq:staticdynamiclink}
\end{equation}
At $u=0$, $\mathcal F_s[0]=0$. Consequently,
\begin{equation}
 -D_u\left.\left( \sqrt{1+\frac{|\nabla u|^2}{(R_s+u)^2}}\mathcal F_s[u]\right)\right|_{u=0}
 =I-\frac{s}{(1+s)\lambda_1(s)}\mathcal I_s:=-\mathcal J_s.
\label{eq:staticdynamicoperators}
\end{equation}
By \eqref{eq:muk}, $-\mathcal J_sY_k^i=\mu_k(s)Y_k^i$.
In particular, the dilation and translation modes have positive
eigenvalues $\mu_0=1$ and $\mu_1=\frac{1}{1+s}$. We remove these modes by
choosing a time-dependent translation and dilation, which leads
to the modulated equation below.

  Set
\begin{equation}
 -\mu_2
 =\frac{s(1+s)-n(1-s)}{(1+s)(n-s)}>0
 \qquad\text{for }s>s_*(n).
 \label{eq:dynamicgap}
\end{equation}
\subsection{Stability for supercritical exponents}
We state the stability theorem. 
\begin{theorem} 
\label{thm:dynamicmain}
Fix $n\geq2$, $s\in(s_*(n),1)$, and $\delta\in(0,1)$. For every $\omega\in(0,-\mu_2)$ there are $\varepsilon,C>0$ with the following property. Suppose that, up to a translation and a dilation, the entire initial boundary is
\[
 \partial E_0=z_0+\rho_0\{(R_s+u_0(\theta))\theta:\theta\in\Sph\},
 \qquad \|u_0\|_{\mathrm h^{2+\delta}}<\varepsilon,
 \quad \rho_0>0,
\]
and $E_0$ is the bounded radial interior of this boundary. Then the classical unconstrained fractional mean curvature flow $V_n=-H_s$ has a finite maximal time $T>0$, is smooth for $0<t<T$, and contracts to a single point $z_\infty$ as $t\uparrow T$.

There exist a time-dependent center $z(\tau)\in\mathbb R^n$,
a scale factor $\rho(\tau)>0$, and an increasing time
reparametrization $t=t(\tau)$, with $t(0)=0$ and
$t(\tau)\to T$ as $\tau\to\infty$, such that
\[
 \partial E_{t(\tau)}
 =z(\tau)+\rho(\tau)
 \bigl\{(R_s+u(\tau,\theta))\theta:
             \theta\in\mathbb S^{n-1}\bigr\}.
\]
The radial height in these moving coordinates satisfies
\begin{equation}\label{eq:dynamicphasebound}
 \|u(\tau)\|_{\mathrm h^{2+\delta}}
 \le C e^{-\omega\tau}
       \|u_0\|_{\mathrm h^{2+\delta}},
 \qquad \tau\ge0.
\end{equation}
Moreover, the normalized profile is smooth for positive times:
for every $k\in\mathbb N$ and $\tau_0>0$, there exists
$C_{k,\tau_0}>0$ such that
\begin{equation}\label{eq:dynamicsmoothbound}
 \|u(\tau)\|_{C^k}
 \le C_{k,\tau_0}e^{-\omega\tau}
       \|u_0\|_{\mathrm h^{2+\delta}},
 \qquad \tau\ge\tau_0.
\end{equation}

The same convergence holds when the flow is rescaled about
its extinction point $z_\infty$ using its extinction time $T$.
More precisely, define
\begin{equation}\label{eq:trueHuisken}
 \widetilde E_q
 =(T-t)^{-1/(1+s)}(E_t-z_\infty),
 \qquad
 q=-\log\frac{T-t}{T}\in[0,\infty),
\end{equation}
so that $q=0$ corresponds to the initial physical time $t=0$.
For every $q\ge0$, the entire rescaled boundary can be written as
\[
 \partial\widetilde E_q
 =\bigl\{(R_s+\widetilde u(q,\theta))\theta:
              \theta\in\mathbb S^{n-1}\bigr\},
\]
where $\widetilde u(q)$ remains small and satisfies
\eqref{eq:dynamicphasebound}--\eqref{eq:dynamicsmoothbound},
with $(u,\tau,\tau_0)$ replaced by
$(\widetilde u,q,q_0)$ and possibly different constants.
Thus the rescaled boundaries converge exponentially to
$\partial B_{R_s}$ in $\mathrm h^{2+\delta}$ and in every
$C^k$ norm for times bounded away from $q=0$. In particular the rescaled limit is unique. The classical second fundamental form satisfies
\begin{equation}
 (T-t)^{1/(1+s)}\|A_{\partial E_t}\|_\infty\leq C,
 \qquad
 (T-t)^{s/(1+s)}\|H_s(E_t)\|_\infty\leq C.
 \label{eq:typeIlocal}
\end{equation}
The constants in the graph estimates do not depend on the initial center or physical scale.
\end{theorem}

Observe that for $0<\delta<\eta$, $C^{j+1,\eta}$ embeds continuously into $\mathrm h^{j+1+\delta}$ since  smoothing approximations converge in the weaker H\"older exponent.  So we get the following immediate corollary. 
\begin{corollary} 
\label{cor:ordinaryHolder}
For any $\eta\in(0,1)$, smallness in $C^{2,\eta}$, up to a similarity, is sufficient for the conclusion of Theorem~\ref{thm:dynamicmain}. One obtains persistence and exponential decay in $\mathrm h^{2+\eta/2}$ and hence in $C^{2,\eta/2}$, as well as convergence in all $C^k$ norms for positive normalized times.
\end{corollary}

We do not claim that these norms are optimal. The extinction point $z_\infty$  and time $T$ are those determined by the solution.
 
We prove Theorem~\ref{thm:dynamicmain} in five steps.

\subsection*{Step 1: Normalization and the radial equation}

For $f$ on $\Sph$, let $\langle f\rangle$ denote its spherical average and define
\begin{equation}
 \mathfrak p f=\bigl(\langle f\rangle,n\langle\theta f\rangle\bigr)
       \in\R\times\R^n,
 \qquad
 \Pi f=\langle f\rangle+n\langle\theta f\rangle\cdot\theta,
 \qquad Q=I-\Pi.
 \label{eq:lowmomentmap}
\end{equation}
Thus $\Pi$ is the projection onto harmonics of degrees $0$ and $1$.

For every positive noninteger $\eta$, define the space of
normalized radial perturbations by
\begin{equation}
 \mathcal X^\eta
 =\left\{u\in\mathrm h^\eta(\Sph):
          \langle u\rangle=0,\ \langle\theta u\rangle=0\right\}
 =\ker\mathfrak p\cap\mathrm h^\eta(\Sph)
 =Q\mathrm h^\eta(\Sph).
 \label{eq:linearslice}
\end{equation}
These conditions remove the degree-$0$ and degree-$1$
components of the radial height. Near the sphere, they
determine a unique choice of translation and dilation,
as shown in Lemma \ref{lem:initialslice}. We refer to these subspaces as the
normalization slice.

\begin{lemma} 
\label{lem:initialslice}
Fix $\eta\in(0,1)$. Every sufficiently small $\mathrm h^{2+\eta}$
radial perturbation of $\partial B_{R_s}$ admits a unique sufficiently
small recentering and rescaling whose radial height $u$ satisfies
$\mathfrak p u=0$.

If $u_0$ denotes the original radial height and $u$ the radial
height after the transformation $E\mapsto\lambda(E-b)$, then
\[
 \|u\|_{\mathrm h^{2+\eta}}+|b|+|\lambda-1|
 \le C\|u_0\|_{\mathrm h^{2+\eta}}.
\]
\end{lemma}
\begin{proof}
Write the boundary of $\lambda(E_{R_s+u_0}-b)$ as a radial
graph about the origin:
\[
 \partial\bigl(\lambda(E_{R_s+u_0}-b)\bigr)
 =\bigl\{(R_s+u(\theta))\theta:   \theta\in\mathbb S^{n-1}\bigr\}
\]
where the new radial height $u$ depends on $(u_0,b,\lambda)$.
We consider the moment  map
\[
 (u_0,b,\lambda)\longmapsto
 \mathfrak p u
 =\bigl(\langle u\rangle,n\langle\theta u\rangle\bigr)
\]
and claim that its  differential with respect to $(b,\lambda)$
at $(u_0,b,\lambda)=(0,0,1)$ is 
\[
 ( b, \lambda)
 \longmapsto (R_s\, \lambda,-  b).
\]
Indeed, let $ b\in\mathbb R^n$ and $ \lambda\in\mathbb R$
be arbitrary parameter variations, and consider
\[
 b_\varepsilon=\varepsilon\, b,
 \qquad
 \lambda_\varepsilon=1+\varepsilon\lambda.
\]
 
By definition $\left|
 \frac{R_s+u_\varepsilon(\theta)}{\lambda_\varepsilon}\theta
 +b_\varepsilon
 \right|^2=R_s^2$ so 
differentiating at $\varepsilon=0$ recalling that 
$u_\varepsilon|_{\varepsilon=0}=0$, yields
\[
 \left.\frac{d}{d\varepsilon}u_\varepsilon(\theta)
 \right|_{\varepsilon=0}
 =R_s \lambda-  b\cdot\theta.
\]

Using the spherical moment identities
$
 \langle\theta\rangle=0,
 \qquad
\langle\theta_i\theta_j\rangle=\frac{\delta_{ij}}n,
$ we obtain
\[
 \left.\frac{d}{d\varepsilon}\mathfrak p u_\varepsilon
 \right|_{\varepsilon=0}
 =
 \left(
 \left\langle R_s \lambda- b\cdot\theta\right\rangle,
 n\left\langle
 \theta(R_s \lambda-b\cdot\theta)
 \right\rangle
 \right)
 =(R_s\lambda,-b).
\]
 
Since $R_s>0$, the differential  is an isomorphism, with inverse
\[
 (\alpha,\beta)\longmapsto
 \left(-\beta,\frac{\alpha}{R_s}\right).
\]
The moment map is $C^1$ in $(u_0,b,\lambda)$.
Indeed, changing variables in its integrals back to the original
radial direction expresses it in terms of smooth functions of
$u_0$, its first derivatives, and the parameters, with denominators
bounded away from zero. The implicit function theorem therefore applies.  

 If $\xi$ is the old radial direction, the new one is the inverse of
\[
 \xi\longmapsto\frac{(R_s+u_0(\xi))\xi-b}{|(R_s+u_0(\xi))\xi-b|}.
\]
For small $u_0,b$ this is a $\mathrm h^{2+\eta}$ diffeomorphism with uniformly controlled inverse. The new radius is $\lambda|(R_s+u_0(\xi))\xi-b|$. The $C^{2+\eta}$ composition and inverse-map estimates give
$\|u\|_{\mathrm h^{2+\eta}}\leq C(\|u_0\|_{\mathrm h^{2+\eta}}+|b|+|\lambda-1|)$; preservation of the little spaces follows by approximation. 
\end{proof}
We now choose the moving center and scale so that
$\mathfrak p u=0$ is preserved. With $r=R_s+u$, set
\begin{equation}
 Z_u=\theta-\frac{\nabla u}{R_s+u}
     =\frac{\nu}{\theta\cdot\nu}.
 \label{eq:radialmobility}
\end{equation}
Let $\tau$ denote time in the moving normalization. Translation
and dilation contribute to its normal velocity through
coefficients $b(u)\in\R^n$ and $c(u)\in\R$, respectively:
\begin{equation}
 V_n=-H_s(r)+c(u)X_u\cdot\nu+b(u)\cdot\nu.
 \label{eq:actualnormalizedV}
\end{equation}
Expressing this speed in radial coordinates and using
\eqref{eq:staticdynamiclink} gives
\begin{equation}
 \begin{aligned}
 u_\tau=-\sqrt{1+\frac{|\nabla u|^2}{(R_s+u)^2}}\,\mathcal F_s[u]+\left(c(u)-\frac1{1+s}\right)r+b(u)\cdot Z_u.
 \end{aligned}
\label{eq:radialmodulated}
\end{equation}
Thus only the last two terms must be determined to preserve the
normalization. Define the $(n+1)\times(n+1)$ matrix $M_u$ by
\begin{equation}
 M_u\binom{c}{b}=\mathfrak p\bigl(cr+b\cdot Z_u\bigr).
\label{eq:modulationmatrix}
\end{equation}
Since $
M_0\binom cb
=\mathfrak p(cR_s+b\cdot\theta)
=
\binom{
\langle cR_s+b\cdot\theta\rangle
}{
n\langle\theta(cR_s+b\cdot\theta)\rangle
}=\binom{
cR_s
}{
b
}$ we get $ M_0=\operatorname{diag}(R_s,I_n)$. Since $M_u\to M_0$ as $u\to0$ in $C^1$, the matrix $M_u$
is invertible for all sufficiently small $u$.

To preserve normalization $\mathfrak p u=0$, we require
$\mathfrak p u_\tau=0$. Applying $\mathfrak p$ to
\eqref{eq:radialmodulated} gives
\[
 M_u
 \binom{c(u)-\frac1{1+s}}{b(u)}
 =
 \mathfrak p\left(
 \sqrt{1+\frac{|\nabla u|^2}{(R_s+u)^2}}\,
 \mathcal F_s[u]
 \right).
\]
Since $M_u$ is invertible for small $u$, this system uniquely
determines $c(u)$ and $b(u)$ for each such profile.  

To write the resulting evolution compactly, we introduce
$\mathsf Q_u$: given a function $\varphi$, it subtracts the unique
combination of $r$ and the components of $Z_u$ having the same
moments as $\varphi$. By the definition of $M_u$, the coefficient
vector of this correction is $M_u^{-1}\mathfrak p\varphi$.
Thus
\begin{equation}
 \begin{aligned}
 \binom{c(u)-\frac1{1+s}}{b(u)}
&=M_u^{-1}\mathfrak p\left(
       \sqrt{1+\frac{|\nabla u|^2}{(R_s+u)^2}}\,\mathcal F_s[u]
       \right),\\
 \mathsf Q_u\varphi
   &=\varphi-( r, Z_u^T)M_u^{-1}\mathfrak p\varphi.
 \end{aligned}
 \label{eq:exactmodulation}
\end{equation}
Indeed, $\mathfrak p\mathsf Q_u=0$ and
$\mathsf Q_u\varphi=\varphi$ when $\mathfrak p\varphi=0$.
Moreover, $\mathsf Q_0=Q$, where $Q$ was introduced in \eqref{eq:lowmomentmap} since 
$\mathsf Q_0\varphi
 =\varphi-\langle\varphi\rangle
-n\langle\theta\varphi\rangle\cdot\theta=Q\varphi$ and $\mathsf Q_u r=0$ since $M_u\binom{1}{0}=\mathfrak p r$.
Hence $\mathsf Q_u$ is a projection onto the fixed space
$\ker\mathfrak p$. Its dependence on $u$ records the dilation
and translation directions of the current profile, represented
by $r$ and the components of $Z_u$.
The closed equation for the normalized height is therefore
\begin{equation}
 u_\tau=-\mathsf Q_u\left(
       \sqrt{1+\frac{|\nabla u|^2}{(R_s+u)^2}}\,\mathcal F_s[u]
       \right):=G(u).
\label{eq:closedshapePDE}
\end{equation}
By the definition of $\mathcal F_s$ in \eqref{eq:staticdynamiclink} and the identity
$\mathsf Q_u r=0$,   equation \eqref{eq:closedshapePDE} is  equivalent to \[u_\tau=-\mathsf Q_u\left(
       \sqrt{1+\frac{|\nabla u|^2}{(R_s+u)^2}}\,H_s(R_s+u)
       \right).\]
 
  \subsection*{Step 2: Local well-posedness and linearization}
\label{step:parabolicsetup}
  
We now verify the hypotheses of the local well-posedness and
linearized stability theorems for quasilinear parabolic equations
in interpolation spaces
\cite[Theorems 1.1 and 1.3]{MWstability}.
For the normalized radial equation, this requires the mapping
estimates established in Section~\ref{sec:radial} and the
sectoriality properties proved below.
Together with the spectral gap for $s>s_*(n)$, these results
yield global existence in normalized time and exponential
convergence to the sphere for sufficiently small normalized
perturbations.
A related quasilinear approach to nonlocal curvature flows
of graphs is developed in \cite[Section 3]{MWgraphs}. 

We first compute the linearization and the modulation coefficients. 
\begin{lemma} \label{lem:quadraticmodulation}
At the sphere, $G(0)=0$, $c(0)=1/(1+s)$, and $b(0)=0$.
For every $v\in\mathcal X^{2+\eta}$,
\begin{equation}
 DG(0)[v]=-\mathcal J_s v,\qquad
 Dc(0)[v]=0,\qquad Db(0)[v]=0.
 \label{eq:modulationlinear}
\end{equation}
Moreover, for $u$ sufficiently small in $C^2$ with
$\mathfrak p u=0$,
\begin{equation}
 \left|c(u)-\frac1{1+s}\right|+|b(u)|
       \le C\|u\|_{C^2}^2.
\label{eq:quadraticmodulation}
\end{equation}
\end{lemma}
\begin{proof}
By Proposition~\ref{prop:cutoffnew}, the curvature is smooth from
a small $C^2$ radial neighborhood into $C^0$. Hence
\eqref{eq:exactmodulation} defines smooth functions $c$ and $b$.
Since $\mathcal F_s[0]=0$, their values at zero follow immediately.
Differentiating the first identity in \eqref{eq:exactmodulation}
and using \eqref{eq:staticdynamicoperators} yields
\[
 \binom{Dc(0)[v]}{Db(0)[v]}
       =M_0^{-1}\mathfrak p(\mathcal J_s v)=0.
\]
The last equality holds because $\mathcal J_s$ preserves harmonic
degrees and $v$ has no components of degrees zero or one (since it is in the kernel of $\mathfrak p$).
Differentiating \eqref{eq:closedshapePDE} similarly gives
$DG(0)[v]=-Q\mathcal J_s v=-\mathcal J_s v$.
For $u\in\ker\mathfrak p$ sufficiently small in $C^2$,
the identities $Dc(0)[u]=Db(0)[u]=0$ and Taylor's formula give
\[
\begin{aligned}
c(u)-\frac1{1+s}
 &=\int_0^1(1-t)\,D^2c(tu)[u,u]\,dt,\qquad 
b(u)
 &=\int_0^1(1-t)\,D^2b(tu)[u,u]\,dt.
\end{aligned}
\]
Since the second derivatives are locally bounded, we obtain \eqref{eq:quadraticmodulation}.
\end{proof}
We also compute the derivative of the projection acting on the
constant curvature term. Differentiating
$\mathsf Q_{\varepsilon v}(R_s+\varepsilon v)=0$ at zero gives
$(D\mathsf Q_0[v])R_s=-Qv=-v$. Since
$\mathcal B_s(R_s)=H_s(R_s)=R_s/(1+s)$, it follows that
\begin{equation}
 (D\mathsf Q_0[v])\mathcal B_s(R_s)=-\frac{v}{1+s}.
 \label{eq:DQcontribution}
\end{equation}

Fix $s,\delta\in(0,1)$ and an integer $j\ge1$. Choose
\[
 \max\{s,\delta\}<\gamma<\min\{1,s+\delta\},
 \qquad \gamma-s<\delta_0<\delta,
\]
and set
\begin{equation}
 X_0=\mathcal X^{j+\gamma-s},
 \qquad X_1=\mathcal X^{j+1+\gamma}.
 \label{eq:Banachcouple}
\end{equation}
The embedding $X_1\hookrightarrow X_0$ is dense and compact.
Define
\begin{align}
 \vartheta&=1-\frac{\gamma-\delta}{1+s},&
 \vartheta_0&=1-\frac{\gamma-\delta_0}{1+s},&
 \xi&=\frac{\delta_0+s-\gamma}{1+s}.
 \label{eq:interpolationindices}
\end{align}
Then $0<\xi<\vartheta_0<\vartheta<1$. Using continuous
interpolation, write
$X_\vartheta=(X_0,X_1)^0_{\vartheta,\infty}$, and similarly
for $X_{\vartheta_0}$ and $X_\xi$. The little-H\"older
interpolation identities  on the normalization slices  (recall that $Q$
is bounded on the corresponding H\"older spaces) yield  
\begin{align}
X_\vartheta&=\mathcal X^{j+1+\delta},&
 X_{\vartheta_0}&=\mathcal X^{j+1+\delta_0},&
 X_\xi&=\mathcal X^{j+\delta_0}.
\label{eq:interpolationidentifications}
\end{align}
For $j=1$, the phase space $X_\vartheta$ is the normalized
subspace of $\mathrm h^{2+\delta}$.

Let $O_{\vartheta_0}\subset X_{\vartheta_0}$ be a sufficiently
small open neighborhood of zero, so that $R_s+u>0$ and $M_u$
is invertible for every $u\in O_{\vartheta_0}$.
For each such $u$, Lemma~\ref{lem:radialmapping} gives
\[
 \mathcal A_s(R_s+u):
 \mathrm h^{j+1+\gamma}\longrightarrow\mathrm h^{j+\gamma-s},
 \qquad
 \mathcal B_s(R_s+u)\in\mathrm h^{j+\delta_0},
\]
with the first map linear and bounded.
Moreover, the factor
$\sqrt{1+|\nabla u|^2/(R_s+u)^2}$ and the components of $Z_u$
belong to $\mathrm h^{j+\delta_0}$.
Since $\delta_0>\gamma-s$, multiplication by these functions
is bounded on both $\mathrm h^{j+\gamma-s}$ and
$\mathrm h^{j+\delta_0}$. Formula~\eqref{eq:exactmodulation}
also gives bounded projections
\[
 \mathsf Q_u:\mathrm h^{j+\gamma-s}\longrightarrow X_0,
 \qquad
 \mathsf Q_u:\mathrm h^{j+\delta_0}\longrightarrow X_\xi.
\]

Using the curvature decomposition \eqref{eq:Hsplit}, we write
the normalized radial equation as
\begin{equation}
 u_\tau=\mathcal L(u)u+f(u),
 \label{eq:abstractdecomposition}
\end{equation}
where, for $v\in X_1$,
\begin{equation}
 \begin{aligned}
 \mathcal L(u)v
 &=-\mathsf Q_u\left(
       \sqrt{1+\frac{|\nabla u|^2}{(R_s+u)^2}}\,
       \mathcal A_s(R_s+u)[v]\right),\\
 f(u)
 &=-\mathsf Q_u\left(
       \sqrt{1+\frac{|\nabla u|^2}{(R_s+u)^2}}\,
       \mathcal B_s(R_s+u)\right).
 \end{aligned}
 \label{definitionl}
\end{equation}
Here $u\in O_{\vartheta_0}$ determines the coefficients,
$\mathcal L(u):X_1\to X_0$ is bounded and linear, and
$f(u)\in X_\xi$. Lemma~\ref{lem:radialmapping} and the smooth
dependence of $M_u^{-1}$ on $u$ yield
\begin{equation}
 (\mathcal L,f)\in
C^\infty\bigl(O_{\vartheta_0},       \mathcal L(X_1,X_0)\times X_\xi\bigr).
\label{eq:abstractsmoothness}
\end{equation}

At $u=0$, the square-root equals $1$ and its first
derivative vanishes. Thus \eqref{eq:DQcontribution} gives
\begin{equation}
 Df(0)[v]
 =-Q\,D\mathcal B_s(R_s)[v]+\frac{v}{1+s},
 \qquad v\in X_{\vartheta_0}.
 \label{eq:linearizedf}
\end{equation}
The last term comes from differentiating $\mathsf Q_u$.
Together with \eqref{eq:Hsplit} and \eqref{eq:As}, this shows
that the full linearization is
$\mathcal L(0)+Df(0)=-\mathcal J_s$ on $X_1$.
 
\begin{lemma}
\label{lem:sectorial}
  For each fixed $u\in O_{\vartheta_0}$, we regard
$\mathcal L(u)$ as an operator on $X_0$ with domain $X_1$.
After shrinking $O_{\vartheta_0}$ if necessary, the associated linear problem
$
 \partial_\tau v=\mathcal L(u)v$, with $ v(0)=v_0\in X_0$
has the solution $v(\tau)=e^{\tau\mathcal L(u)}v_0\in X_1$
for $\tau>0$, with $v(\tau)\to v_0$ in $X_0$ as $\tau\downarrow0$.
More precisely, $\mathcal L(u)$ generates a strongly continuous
analytic semigroup on $X_0$ with domain $X_1$.
 
  The linearization of \eqref{eq:abstractdecomposition} at zero is
\begin{equation}
 \mathcal L(0)+Df(0)=-\mathcal J_s|_{X_1},\qquad
 \sup\{\operatorname{Re}\lambda:\lambda\in\sigma(-\mathcal J_s|_{X_0})\}
       =\mu_2.
 \label{eq:spectrumBanach}
\end{equation}
\end{lemma}
\begin{proof}
Lemma~\ref{lem:spectralresolvent}, applied with
$b=j+\gamma-s$, shows that $-\mathcal I_s$ generates a
strongly continuous analytic semigroup on
$\mathrm h^{j+\gamma-s}$ with domain
$\mathrm h^{j+1+\gamma}$.
Since $\mathcal I_s$ preserves spherical harmonic degrees,
its semigroup preserves the normalization slice.
Consequently, the same generation property holds for
$-R_s^{-1-s}\mathcal I_s$ on $X_0$, with domain $X_1$.

We next show that, on $X_1$,
\[
 \mathcal L(0)
 =-R_s^{-1-s}\mathcal I_s
   +\frac{s}{1+s}I+D\mathcal B_s(R_s),
\]
where $D\mathcal B_s(R_s)$ extends to a bounded operator
on $X_0$. 
To identify $\mathcal L(0)$, let
\[
 K_s u(\theta)=\int_\Sph \frac{u(\eta)}{|\theta-\eta|^{n+s-2}} \dd\eta,
 \qquad  K_s1=sH_s(\mathbb S^{n-1})=\lambda_1(s).
\]
This is an integrable rotation-invariant kernel, and $K_s$ is bounded on every noninteger H\"older space: commute rotation derivatives through the integral and use its finite $L^1$ norm for the H\"older differences. Differentiation of \eqref{eq:radiallowerterm} at $r=R_s$ gives the explicit bounded operator
\begin{equation}
 D\mathcal B_s(R_s)[u]=\frac{R_s^{-1-s}}{2s}
 \bigl((2-n-s)\lambda_1(s) u+(n-4-s)K_su\bigr).
 \label{eq:DBsphere}
\end{equation}
By \eqref{eq:Hsplit} and \eqref{eq:As},
\[
 \mathcal A_s(R_s)=\mathcal J_s+\frac1{1+s}I-D\mathcal B_s(R_s)
       =R_s^{-1-s}\mathcal I_s+\text{a bounded operator on }X_0.
\]
All these operators preserve the spherical harmonic subspaces. At $u=0$, $\mathsf Q_0=Q$ and consequently, recalling \eqref{definitionl}, \[\mathcal L(0)=-\mathcal A_s(R_s) \qquad \text{on the slice.}
\]  
Thus $\mathcal L(0)$ differs from
$-R_s^{-1-s}\mathcal I_s$ by an operator bounded on $X_0$.
The bounded perturbation theorem therefore implies that
$\mathcal L(0)$ also generates a strongly continuous analytic semigroup on $X_0$, with the same domain $X_1$. Finally, the continuity established in
\eqref{eq:abstractsmoothness} gives
\[
 \|\mathcal L(u)-\mathcal L(0)\|_{\mathcal L(X_1,X_0)}
 \longrightarrow0
 \qquad\text{as }u\to0\text{ in }X_{\vartheta_0}.
\]
The small-perturbation theorem then yields the same generation
property for $\mathcal L(u)$ after shrinking
$O_{\vartheta_0}$.
In this last step, the perturbation is small as an operator
from $X_1$ to $X_0$; boundedness on $X_0$ is not required.
 
The semigroup estimates give, uniformly for small $u$, $
\|e^{\tau\mathcal L(u)}v_0\|_{X_0}
\le C e^{\kappa\tau}\|v_0\|_{X_0}$, but need not imply decay. Under the above parabolic
hypotheses, exponential stability follows from the linearized
stability theorem provided that the full linearization
$\mathcal L(0)+Df(0)=-\mathcal J_s$ has strictly negative
spectral bound.  
Combining \eqref{eq:linearizedf} with \eqref{eq:Hsplit} gives
\[
 \mathcal L(0)v+Df(0)[v]
   =-Q\,DH_s(R_s)[v]+\frac{v}{1+s}=-\mathcal J_s v,
 \qquad v\in X_1.
\]
This agrees with \eqref{eq:modulationlinear}. Its resolvent is compact. If $(-\mathcal J_s-\lambda)v=0$, applying $P_k$ gives
$(\mu_k-\lambda)P_kv=0$ for each $k\ge2$.
Completeness of the spherical harmonics shows that a nonzero
eigenfunction can occur only when $\lambda=\mu_k$ for some $k$.
Conversely, every spherical harmonic of degree $k\ge2$ lies
in the domain. The spectrum on the slice is consequently exactly $\{\mu_k:k\geq2\}$. Since $\lambda_k$ is strictly increasing, we conclude  \eqref{eq:spectrumBanach}.
\end{proof}

  \subsection*{Step 3: Global existence, decay, and regularization}
\label{step:globalexistence}
  
Global existence for the normalized radial equation
\eqref{eq:abstractdecomposition} and estimate
\eqref{eq:dynamicphasebound} follow from the result below,
quoted from \cite[Theorems 1.1 and 1.3, and Remark 1.2(iii)]{MWstability}.
We state its specialization to the zero equilibrium, with the
generator sign convention of \eqref{eq:abstractdecomposition}.
Its proof uses parabolic evolution-operator estimates in
interpolation spaces and a continuation argument.  

\begin{theorem}
\label{thm:abstractinput}
Let $X_1$ be densely embedded in $X_0$, use admissible
interpolants $X_\theta$, and suppose
$0<\xi\leq\vartheta_0<\vartheta<1$.
On an open neighborhood of zero in $X_{\vartheta_0}$, assume
that $\mathcal L$ is locally Lipschitz into the analytic generators
on $X_0$ with common domain $X_1$, and that $f$ is locally
Lipschitz into $X_\xi$.
Suppose that $f(0)=0$ and that $f$, regarded as a map into
$X_0$, is Fr\'echet differentiable at zero, with
$Df(0)\in\mathcal L(X_{\vartheta_0},X_0)$. Assume also that
\[
 \sigma\bigl(\mathcal L(0)+Df(0)\bigr)
      \subset\{\operatorname{Re}z\leq-\omega_0\},\qquad\omega_0>0.
\]
Then, for every $0<\omega<\omega_0$, sufficiently small initial data in $X_\vartheta$ generate unique global solutions of $u'=\mathcal L(u)u+f(u)$ and
\[
 \|u(\tau)\|_{X_\vartheta}
       \leq C e^{-\omega\tau}\|u(0)\|_{X_\vartheta}.
\]
The solutions are continuous into $X_\vartheta$ at zero, and belong to $C((0,\infty),X_1)\cap C^1((0,\infty),X_0)$. Smooth coefficient maps give smooth dependence at positive time on data and on finite-dimensional parameters that preserve this fixed Banach couple and its generator domain.  
\end{theorem}
 
The assumptions of the theorem are satisfied in our setting; density and compactness follow from \eqref{eq:Banachcouple}; the strict interpolation inequalities are \eqref{eq:interpolationindices}; coefficient smoothness and the positive regularity of $f$ are \eqref{eq:abstractsmoothness}; sectoriality is Lemma~\ref{lem:sectorial}; the equilibrium and its derivative are Lemma~\ref{lem:quadraticmodulation}; and the spectral bound is the explicit number \eqref{eq:dynamicgap}. 

Theorem~\ref{thm:abstractinput} therefore yields a global
solution of the normalized radial equation \eqref{eq:abstractdecomposition} satisfying
\eqref{eq:dynamicphasebound}.  
For sufficiently small initial data, this estimate keeps
the solution within the neighborhood where the equation
is defined. We also record the following continuation
criterion, which does not require $s>s_*(n)$.

\begin{lemma} 
\label{lem:continuation}
Fix $s\in(0,1)$ and use the spaces and interpolation
parameters defined in
\eqref{eq:Banachcouple}--\eqref{eq:interpolationindices}.
A maximal solution of the normalized radial equation
\eqref{eq:abstractdecomposition} exists for all normalized times
$\tau\ge0$ provided that, throughout its interval of existence,
its $X_\vartheta$ norm remains bounded and its distance
from $\partial O_{\vartheta_0}$, measured in
$X_{\vartheta_0}$, is bounded below by a positive constant.
In particular, any solution remaining in a sufficiently
small fixed ball about zero in $X_\vartheta$ is global.
\end{lemma}

\begin{proof}
The embedding $X_1\hookrightarrow X_0$ is compact, since
the sphere is compact and the regularity indices differ
by $1+s>0$.
As $\vartheta_0<\vartheta<1$, the compact-embedding
alternative in \cite[Theorem~1.1(v)]{MWstability} applies
with intermediate exponent $\vartheta$.
A finite maximal time would therefore force either
the $X_\vartheta$ norm to become unbounded or the
$X_{\vartheta_0}$ distance to $\partial O_{\vartheta_0}$
to approach zero. Both possibilities are excluded.

Finally, the continuous embedding
$X_\vartheta\hookrightarrow X_{\vartheta_0}$ ensures that
a sufficiently small ball about zero in $X_\vartheta$
lies at a positive $X_{\vartheta_0}$ distance from
$\partial O_{\vartheta_0}$. A solution confined to this
ball satisfies both conditions above.
\end{proof}

 We next establish positive-time regularity and decay in higher norms. 

\begin{lemma}
\label{lem:positivetimesmoothing}
Let $u$ be the global solution of the normalized radial
equation obtained for sufficiently small initial data.
Then
\[
 u\in C^\infty\bigl((0,\infty)\times\mathbb S^{n-1}\bigr),
\]
and, for every $k\in\mathbb N$ and $\tau_0>0$,
estimate \eqref{eq:dynamicsmoothbound} holds.
The constants depend on $k$, $\tau_0$, the fixed parameters,
and the chosen neighborhood, but not on the initial datum.
The smallness condition on the initial datum is independent
of $k$ and $\tau_0$.
\end{lemma}

\begin{proof}
We first establish regularization and a uniform estimate
at a fixed positive time. Throughout the proof, take $j=1$,
so that $X_\vartheta=\mathcal X^{2+\delta}$.

The normalized equation and the constraint
$\mathfrak p u=0$ are invariant under rotations.
For a skew-symmetric matrix $K$, define
$u_K(\tau,\theta)=u(\tau,e^{\tau K}\theta).$
Then $u_K(0)=u(0)$ and
$
 \partial_\tau u_K=G(u_K)+(K\theta)\cdot\nabla u_K.$ 
The added operator preserves the constraint and is bounded
from $X_1$ to $X_0$, with norm $O(|K|)$.
Hence, for small $K$, the modified equation satisfies
the same parabolic hypotheses.
The smooth-dependence result in
\cite[Remark~1.2(iii)]{MWstability}, applied with $K$
as an additional parameter, gives smooth dependence
on $K$, the initial datum, and positive time.

At a fixed $\tau>0$, differentiation with respect to $K$
produces spatial derivatives. For example,
\[
 \left.\frac{d}{d\varepsilon}
 u_{\varepsilon K}(\tau,\theta)\right|_{\varepsilon=0}
 =\tau(K\theta)\cdot\nabla u(\tau,\theta).
\]
Higher parameter derivatives similarly control higher
spatial derivatives: rotations provide local coordinates
on the sphere, since their infinitesimal generators span
each tangent space.
Together with smooth dependence on positive time,
this proves smoothness in space and time.
This is the parameter argument of
\cite[Section~3.4]{MWgraphs}, with rotations replacing
translations.

To obtain a quantitative estimate, let $S(h)w$ denote
the solution at time $h$ with initial datum $w$.
The uniform coefficient bounds in
Lemma~\ref{lem:radialmapping}, together with the modulation
formulas, allow the local solution construction on a fixed
small data ball and a common time interval.
Differentiating its fixed-point equation with respect
to the rotation parameters and the initial datum always
uses the same uniformly invertible linear operator.
Induction therefore bounds derivatives of every fixed
order on this same ball; only the constants depend
on the order.

For initial datum identically $0$, the solution is identically $0$
for every $K$. The mean value theorem in the initial datum,
combined with the preceding parameter estimates, thus gives
\[
 \|S(h)w\|_{C^k}
 \le C_{k,h}\|w\|_{X_\vartheta}
\]
for every positive $h$ in the common local existence
interval and every $w$ in the fixed small ball.

Now fix $\tau_0>0$ and choose $0<h<\tau_0$ within that
interval. By \eqref{eq:dynamicphasebound}, sufficiently
small initial data ensure that the entire trajectory
remains in the ball where the preceding estimate applies.
Using uniqueness and the autonomous character of the
equation, we obtain, for $\tau\ge\tau_0$,
\[
 \begin{aligned}
 \|u(\tau)\|_{C^k}
  =\|S(h)u(\tau-h)\|_{C^k} \le C_{k,h}\|u(\tau-h)\|_{X_\vartheta} \le C_{k,h}C e^{\omega h}e^{-\omega\tau}
       \|u_0\|_{X_\vartheta}.
 \end{aligned}
\]
 
\end{proof}

\subsection*{Step 4: Reconstruction of the physical flow and extinction}

The radial equation is written in normalized time $\tau$; it
remains to check that it describes the original \emph{unforced}
fractional mean curvature flow and that its physical extinction time
is finite.  Let $u$ be the global solution of the normalized radial equation, with initial data selected as in Lemma~\ref{lem:initialslice}. For chosen physical scale $\rho_0>0$ and center $z_0$, solve the scalar/vector equations (where the derivatives are w.r.t. $\tau$):
\begin{equation}
 \rho'  =-c(u)\rho,\qquad z'=-\rho b(u),\qquad
 t'=\rho^{1+s},\qquad
 \rho(0)=\rho_0,\ z(0)=z_0,\ t(0)=0.
 \label{eq:reconstructionODE}
\end{equation}
  Define
$ E_{t(\tau)}=z(\tau)+\rho(\tau)E_{R_s+u(\tau)}$ so that  $\partial E_{t(\tau)}$ is parametrized by
$z(\tau)+\rho(\tau)X_u(\tau,\theta)$, where
$X_u(\tau,\theta)=(R_s+u(\tau,\theta))\theta$.
Differentiating this parametrization with respect to $\tau$
at fixed $\theta$, and taking its normal component, gives
\[
\begin{aligned}
\partial_\tau\bigl(z+\rho X_u\bigr)\cdot\nu
&=z'\cdot\nu+\rho'X_u\cdot\nu
  +\rho\,\partial_\tau X_u\cdot\nu\\
&=z'\cdot\nu+\rho'X_u\cdot\nu
  +\rho\bigl(-H_s+cX_u\cdot\nu+b\cdot\nu\bigr) =-\rho H_s,
\end{aligned}
\]
where the last equality follows from
$z'=-\rho b$ and $\rho'=-\rho c$.
Here $\nu$ is the common outer unit normal to the normalized
and physical boundaries at corresponding points, and $H_s$
is evaluated on the normalized boundary.
Dividing by $t'=\rho^{1+s}$ gives $-\rho^{-s}H_s$, which is exactly minus the fractional curvature of the reconstructed set. Thus  $E_{t(\tau)}$ solves $V_n=-H_s$. 

For small data, $c(u)\geq 1/(2(1+s))$ and $H_s(E_r)$ is uniformly positive and bounded, since curvature is continuous in a $C^2$ neighborhood of the sphere. Hence $\rho$ stays positive, tends to zero exponentially, and both $T$ and $z_\infty$, defined as 
\begin{equation}
T:=\int_0^\infty\rho(\tau)^{1+s}\dd\tau<\infty,
 \qquad
 z_\infty:=z_0-\int_0^\infty\rho(\tau)b(u(\tau))\dd\tau
 \label{eq:extinctionparameters}
\end{equation}
are well defined. The  flow is classical on every compact subinterval of $[0,T)$ and smooth at positive times. Local uniqueness, given by classical short-time theory  \cite{JLMshort}, identifies it with the original classical solution. Its boundary collapses to $z_\infty$, so it cannot be continued through $T$ as a classical embedded hypersurface. Thus $T$ is its maximal classical time.  
 
Write $e_0=\|u(0)\|_{X_\vartheta}$. Lemma~\ref{lem:quadraticmodulation} and the phase estimate \eqref{eq:dynamicphasebound} give
\begin{equation}
 \left|c(u(\tau))-\frac1{1+s}\right|+|b(u(\tau))|
      \leq C e_0^2 e^{-2\omega\tau}.
 \label{eq:modulationdecay}
\end{equation}
In particular
\begin{equation}
\rho_\infty=\rho_0\exp\!\left(-\int_0^\infty\left(c(u)-\frac1{1+s}\right)\dd\tau\right)<\infty
 \label{eq:rhoinfinity}
\end{equation}
 and direct integration yields
\begin{align}
\rho(\tau)&=\rho_\infty e^{-\tau/(1+s)}
\bigl(1+O(e_0^2e^{-2\omega\tau})\bigr),
\label{eq:rhoasymptotic}\\
 T-t(\tau)&=\rho_\infty^{1+s}e^{-\tau}
\bigl(1+O(e_0^2e^{-2\omega\tau})\bigr),
\label{eq:timeasymptotic}\\
 \frac{|z(\tau)-z_\infty|}{\rho(\tau)}
       &\leq C e_0^2e^{-2\omega\tau},
\label{eq:centerasymptotic}\\
 \frac{\rho(\tau)}{(T-t(\tau))^{1/(1+s)}}
       &=1+O(e_0^2e^{-2\omega\tau}).
\label{eq:scaleasymptotic}
\end{align}
For example, \eqref{eq:centerasymptotic} follows by integrating
$\rho(\tau')|b(u(\tau'))|\leq C\rho_\infty e^{-(\frac1{1+s}+2\omega)\tau'}e_0^2$
from $\tau$ to infinity and dividing by \eqref{eq:rhoasymptotic}. In \eqref{eq:timeasymptotic}, the leading term of $\rho^{1+s}$ is $\rho_\infty^{1+s}e^{-\tau}$, whose integral from $\tau$ to infinity has the same value.

For initial data already in the slice, these formulas also give
\begin{equation}
 T=\rho_0^{1+s}\bigl(1+O(e_0^2)\bigr),
 \qquad |z_\infty-z_0|\leq C\rho_0e_0^2.
 \label{eq:quadraticpolechange}
\end{equation}
The quadratic estimates above are relative to the initial
center and scale selected in Lemma~\ref{lem:initialslice}.
For a general initial surface close to a sphere, achieving
this normalization may require a translation and a change
of scale of first order in the original perturbation.
Only the subsequent corrections to the extinction parameters
are quadratic in the normalized perturbation.

Let $q$ be as in \eqref{eq:trueHuisken}, consistently with
the fixed-pole rescaling $q=-\log((T-t)/T)$; then $q(0)=0$.  By \eqref{eq:timeasymptotic},
\[
 q=\tau+\log\frac{T}{\rho_\infty^{1+s}}
                   +O(e_0^2e^{-2\omega\tau}).
\]
Thus the two normalized times differ by a bounded constant and a decaying error. Moreover
\begin{equation}
 \widetilde E_q
 =\frac{z(\tau)-z_\infty}{(T-t)^{1/(1+s)}}
     +\frac{\rho(\tau)}{(T-t)^{1/(1+s)}}E_{R_s+u(\tau)}.
 \label{eq:gaugecomparison}
\end{equation}
The translation and extra dilation in this formula are $O(e_0^2e^{-2\omega\tau})$. Writing the transformed boundaries again as radial graphs,
as in Lemma~\ref{lem:initialslice}, and using the same composition
estimates transfers \eqref{eq:dynamicphasebound} and
\eqref{eq:dynamicsmoothbound} to the fixed-pole graphs. In this fixed-pole normalization the speed is precisely $-\mathfrak D_s$.
\subsection*{Step 5: Conclusion}
\begin{proof}[Completion of the proof of Theorem~\ref{thm:dynamicmain}]
Select the initial slice by Lemma~\ref{lem:initialslice}. Theorem~\ref{thm:abstractinput}, with all hypotheses verified in \hyperref[step:parabolicsetup]{Step~2}, gives the global solution of the normalized radial equation and the phase estimate. Lemma~\ref{lem:positivetimesmoothing} gives every higher norm on positive-time strips. Equations \eqref{eq:reconstructionODE}--\eqref{eq:extinctionparameters} give the original flow, its finite maximal time, and its extinction point. Equations \eqref{eq:rhoasymptotic}--\eqref{eq:gaugecomparison} prove convergence and persistence in the true Huisken normalization. Finally $A$ scales with degree $-1$, $H_s$ with degree $-s$, and the normalized $C^2$ geometry is uniformly controlled; this proves \eqref{eq:typeIlocal}.
\end{proof}

\subsection{Instability for subcritical exponents}
The local well-posedness and smooth-dependence results
established in \hyperref[step:parabolicsetup]{Step~2} hold for every
$s\in(0,1)$. The condition $s>s_*(n)$ was used only to obtain
exponential stability. For $s<s_*(n)$, we instead apply the
principle of linearized instability to the same normalized
radial equation, which is contained in \cite[Theorem~1.4]{MWstability}.
\begin{theorem}
\label{prop:belowthreshold}
For $0<s<s_*(n)$, the sphere is nonlinearly unstable modulo translations and dilations in the same little-H\"older phase spaces. 
\end{theorem}
No assertion is made at $s=s_*(n)$, where the degree-two eigenvalue vanishes and the
linearized instability theorem used here does not apply.
\begin{proof}
On the normalization slice, the eigenvalue of
$-\mathcal J_s$ associated with spherical harmonics
of degree two is
$\mu_2(s)=\frac{n-(n+1)s-s^2}{(n-s)(1+s)}>0$
 for $0<s<s_*(n)$.
Thus the linearized equation has exponentially growing
degree-two modes. Moreover, $\mu_k(s)\to-\infty$ as
$k\to\infty$, so the unstable spectrum is nonempty
and consists of finitely many eigenvalues.
The smoothness established in \eqref{eq:abstractsmoothness}
implies  the differentiability assumptions of
\cite[Theorem~1.4]{MWstability}. In particular, the embedding
$X_\xi\hookrightarrow X_0$ gives
\[
 Df(0)\in\mathcal L(X_{\vartheta_0},X_0),
 \qquad 0<\vartheta_0<1,
\]
so the derivative of the lower-order term is bounded from
a strict intermediate space into the base space.
The theorem therefore yields arbitrarily small initial
perturbations whose solutions of the normalized radial
equation leave a fixed neighborhood of zero.

The initial perturbations may also be chosen smooth.
Indeed, smooth functions are dense in the phase space
within the slice, and continuous dependence at a fixed
exit time ensures that sufficiently close approximations
still leave a slightly smaller neighborhood.
Since the normalization removes translations and dilations,
this proves nonlinear instability modulo these symmetries.

\end{proof}

\section{Lyapunov--Schmidt reduction}\label{sec:LS}

The preceding section establishes the stability or instability of the sphere
away from the critical exponent. We now study the shrinker equation
near $\sigma=s_*(n)$. For the reduction in this section, write
$v=u/R_s$ and normalize the stationary defect by
\begin{equation}\label{eq:staticF}
 \mathscr F(s,v):=\frac{1+s}{R_s}\mathcal F_s[R_sv]
 =\frac{H_s(1+v)}{H_s(\mathbb S^{n-1})}
 -\frac{1+v}{\sqrt{1+|\nabla v|^2/(1+v)^2}}.
\end{equation}
Here $H_s(1+v)(\theta)=H_s(E_{1+v},(1+v(\theta))\theta)$,
in the radial notation of Section~\ref{sec:radial}. The equality
follows from the homogeneity of $H_s$ and
$R_s^{1+s}=(1+s)H_s(\mathbb S^{n-1})$.
Thus $\mathscr F(s,v)=0$ is equivalent to
$\mathcal F_s[R_sv]=0$, and $\mathscr F(s,0)=0$.

By the chain rule and \eqref{eq:As}, its linearization is
\begin{equation}\label{eq:physicalstaticlinear}
 \mathscr L_s:=D_v\mathscr F(s,0)
 =(1+s)\mathcal J_s
 =\frac{s}{\lambda_1(s)}\mathcal I_s-(1+s)I.
\end{equation}
Consequently, on spherical harmonics of degree $k$,
\begin{equation}\label{eq:staticlinear}
 \mathscr L_sY_k^i=\ell_k(s)Y_k^i,\qquad
 \ell_k(s)=s\frac{\lambda_k(s)}{\lambda_1(s)}-(1+s)
 =-(1+s)\mu_k(s).
\end{equation}
The constant and translation modes are not kernel directions:
\(\ell_0=-(1+s)\) and \(\ell_1=-1\). Moreover,
\begin{equation}\label{eq:transversality}
 \ell_2(s)=\frac{s^2+(n+1)s-n}{n-s},\qquad
 \ell_2(\sigma)=0,\qquad
 d_n:=\ell_2'(\sigma)=\frac{2n^2}{(n-\sigma)^2}-1>0.
\end{equation}
The strict increase of \(\lambda_k(s)\) in \(k\) shows that
\(\ker\mathscr L_\sigma=\mathcal H_2\).

The kernel has dimension $n(n+1)/2-1>1$, so the
simple-eigenvalue bifurcation theorem of Crandall and
Rabinowitz \cite[Theorem~1.7]{CR} cannot be applied directly
to the unrestricted equation.
We therefore first perform a Lyapunov--Schmidt reduction
onto $\mathcal H_2$. Restricting the reduced equation to
suitable one-dimensional symmetry-fixed subspaces and
factoring out the trivial branch yields scalar equations.
The transversality condition \eqref{eq:transversality}
then allows us to construct the bifurcating branches
by the scalar implicit function theorem.

Write \(P_2\) for the \(L^2\)-orthogonal projection onto
\(\mathcal H_2\), and \(Q_2=I-P_2\). These projections act
boundedly on the little-H\"older spaces below.

Choose \(\eta\in(\sigma,1)\) and a compact interval
\(I\Subset(0,\eta)\) containing \(\sigma\) in its interior. For
integers \(r\ge2\), set
\[
 X_r=\mathrm h^{r+\eta}(\mathbb S^{n-1}),\qquad
 Y_{r,s}=\mathrm h^{r-1+\eta-s}(\mathbb S^{n-1}).
\]
The loss of \(1+s\) derivatives is explicit in this notation.

\begin{lemma} 
\label{lem:ellipticscale}
For each fixed \(s\in I\), the map
\(\mathscr F(s,\cdot):X_r\to Y_{r,s}\) is smooth near zero.
The estimates are locally uniform for \(s\in I\). Moreover,
\begin{equation}\label{eq:complementinverse}
 \mathscr B_s:=(Q_2\mathscr L_s|_{Q_2X_r})^{-1}:
 Q_2Y_{r,s}\longrightarrow Q_2X_r
\end{equation}
exists and is uniformly bounded for \(s\) in a sufficiently small
closed interval about \(\sigma\). 

Let
\(\mathscr N_s(v)=\mathscr F(s,v)-\mathscr L_sv\), then
\begin{equation}\label{eq:Nlipschitz}
 \|\mathscr N_s(v)-\mathscr N_s(z)\|_{Y_{2,s}}
 \le C(\|v\|_{X_2}+\|z\|_{X_2})\|v-z\|_{X_2}
\end{equation}
for sufficiently small \(v,z\).
\end{lemma}
\begin{proof}
Smoothness follows from the radial estimates of
Section~\ref{sec:radial}, with input H\"older exponent $\eta>s$;
the support term has lower order. Taylor's formula gives
\eqref{eq:Nlipschitz}, uniformly for $s$ in the indicated compact interval.

The inverse, its uniformity and its Schauder bounds follow
from Lemma~\ref{lem:spectralresolvent}, with
$b=r-1+\eta-s$.
\end{proof}

We identify \(\mathcal H_2\) with the space
\(\operatorname{Sym}_0(n)\) of trace-free symmetric matrices by setting
\begin{equation}\label{eq:rangeparameters}
q_A(\theta)=\theta^TA\theta,\qquad A\in\operatorname{Sym}_0(n).
\end{equation}

We seek solutions in the form $v=q_A+w$, where
$q_A\in\mathcal H_2$ and $P_2w=0$. The range equation will
determine $w=w(s,A)$, and the full equation becomes
\begin{equation}\label{deco}\mathscr F(s,v)=0\text{ is equivalent to  }\begin{cases}P_2\mathscr F(s,q_A+w(s,A))=0\\ Q_2\mathscr F(s,q_A+w(s,A))=0.\end{cases}
\end{equation}
The next lemma solves the second equation in \eqref{deco}
and establishes the joint regularity of the range correction. \begin{lemma} \label{lem:range}
On a neighborhood of $(\sigma,0)$ there is a unique small
$w(s,A)\in Q_2X_2$, $A\in\operatorname{Sym}_0(n)$, satisfying
\begin{equation}\label{eq:rangeeq}
 Q_2\mathscr F(s,q_A+w(s,A))=0.
\end{equation}
The map is jointly $C^\infty$ into every $X_r$, and
\begin{equation}\label{wsecondo}
 w(s,0)=D_Aw(s,0)=0,\qquad \|w(s,A)\|_{X_r}\le C_r|A|^2.
\end{equation}
The neighborhood is independent of $r$; only $C_r$ depends on $r$.
Moreover
\begin{equation}\label{eq:rangeorbit}
 U_Rw(s,A)=w(s,RAR^T),\qquad U_Rf(\theta)=f(R^{-1}\theta).
\end{equation}
In particular $w(s,A)$ is even.
\end{lemma}
\begin{proof}   For fixed \(s\), the range equation \eqref{eq:rangeeq} is equivalent to
\begin{equation}\label{eq:rangecontraction}
 w=-\mathscr B_sQ_2\mathscr N_s(q_A+w).
\end{equation}
 
By Lemma~\ref{lem:ellipticscale}, the right-hand side of
\eqref{eq:rangecontraction} maps the ball
$\|w\|_{X_2}\le C|A|^2$ into itself and is a contraction with
constant less than $1/2$, uniformly in $s$ after shrinking
$I$ and the low-norm neighborhood. The same Lipschitz bound
gives uniqueness among all sufficiently small range corrections.
For fixed $s$, the implicit function theorem gives smooth
dependence on $A$, with uniform derivative bounds.
 
Equivariance and uniqueness give \eqref{eq:rangeorbit}. Taking $R=-I$ proves evenness. For a
skew-symmetric matrix $K$, differentiation of \eqref{eq:rangeorbit} along $e^{tK}$ gives
\[
 \left.\frac{\dd}{\dd t}\right|_0 U_{e^{tK}}w(s,A)
 =D_Aw(s,A)[KA-AK].
\]
Iterated rotation fields span all spatial derivatives on the sphere.
Thus $w$ is spatially smooth. The bounds remain $O_r(|A|^2)$:
$D_Aw=O(|A|)$ multiplies an $O(|A|)$ argument, while terms with two
or more $A$-derivatives have at least two such arguments. The same
reasoning gives uniform mixed spatial and $A$-derivative bounds.
Up to this point, no continuity in $s$ has been used.

To prove joint regularity, set
$T_{s,v}=Q_2D_v\mathscr F(s,v)Q_2$, regarded as a map
$Q_2X_2\to Q_2Y_{2,s}$.
Since $D\mathscr N_s(0)=0$, the low-norm neighborhood can
be chosen so that
\[
 \|\mathscr B_sQ_2D\mathscr N_s(v)Q_2\|_{\mathcal L(Q_2X_2,Q_2X_2)}
 \le \frac12
\]
uniformly in $s$. The Neumann series therefore gives
\[
 T_{s,v}^{-1}
 =\bigl(I+\mathscr B_sQ_2D\mathscr N_s(v)Q_2\bigr)^{-1}
   \mathscr B_s.
\]
The same bound holds for averages along segments in this
neighborhood. For smooth $v$ and smooth $f\in Q_2Y_{2,s}$,
equivariance gives
\[
 U_R(T_{s,v}^{-1}f)=T_{s,U_Rv}^{-1}U_Rf.
\]
 
The right-hand side depends smoothly on $R$ in the low pair,
because $v,f$ are smooth and inversion is smooth on bounded
isomorphisms. Consequently the rotation orbit of
$z=T_{s,v}^{-1}f$ is smooth in $X_2$, which justifies
differentiating before any higher regularity of $z$ is assumed.
For a rotation generator $Z$, the differentiated equation is
\[
 T_{s,v}(Zz)=Zf-Q_2D_v^2\mathscr F(s,v)[Zv,z].
\]
Iteration proves that smooth data give a smooth solution. More
quantitatively, for every integer $m\ge0$,
\[
 \|T_{s,v}^{-1}f\|_{X_{2+m}}
 \le C_m\|f\|_{Y_{2+m,s}},
\]
when the finitely many higher norms of $v$ occurring in the induction
are bounded. Indeed after $m$ rotations the right-hand side is in
$Y_{2,s}$ by the fixed-order multilinear estimates, using the previously
controlled rotation derivatives of $z$. The same argument holds for
averages of these operators along segments of uniformly smooth small
profiles.  

 Proposition~\ref{prop:cutoffnew}
gives joint smoothness with the fixed loss $C^{k+2}\to C^k$.
The already uniformly smooth profiles make this sufficient for every
finite target norm. In particular, for
$v_s=q_A+w(s,A)$ and $\Delta_hw=w(s+h,A)-w(s,A)$, subtraction gives
\begin{align}
 T_{s,h}\Delta_hw
 &=-Q_2[\mathscr F(s+h,v_{s+h})-\mathscr F(s,v_{s+h})],
 \label{eq:rangedifference}\\
 T_{s,h}&=Q_2\int_0^1D_v\mathscr F(s,v_s+t\Delta_hw)\dd t\,Q_2.
 \notag
\end{align}
The order in the averaged linearization is $s$. The right-hand side
is $O(|h|)$ in every prescribed finite norm, uniformly by the cutoff
estimate with logarithms and the spatial bounds above. The inverse
estimates give $\|\Delta_hw\|_{X_r}\le C_r|h|$ for every $r$.
Dividing by $h$ and comparing the resulting equations in the low
pair yields
\begin{equation}\label{eq:rangeparameterderivative}
 T_{s,v_s}\partial_sw=-Q_2\partial_s\mathscr F(s,v_s).
\end{equation}
Convergence in higher norms follows by interpolation with the uniform
bounds for the difference quotients in arbitrarily high norms.
Each subsequent mixed derivative satisfies an equation with the same
invertible $T_{s,v_s}$ and a right-hand side containing already
constructed derivatives and derivatives of $\mathscr F$ on smooth
directions. For example
\[
 T_{s,v_s}w_{ss}=-Q_2\bigl(\mathscr F_{ss}
 +2D_v\mathscr F_s[w_s]+D_v^2\mathscr F[w_s,w_s]\bigr).
\]
Repeating the difference-quotient argument justifies these identities
and joint continuity inductively. This proves joint smoothness on the
same low-norm neighborhood.  
\end{proof}


By Lemma~\ref{lem:range},
$\mathscr F(s,q_A+w(s,A))$ belongs to $\mathcal H_2$.
The identification \eqref{eq:rangeparameters} therefore defines
a unique trace-free symmetric matrix $g(s,A)$ by
\begin{equation}\label{eq:reduceddefinition}
 q_{g(s,A)}=P_2\mathscr F(s,q_A+w(s,A)).
\end{equation}
Its zeros correspond exactly to the sufficiently small
solutions of the full shrinker equation.
The map $g$ is jointly smooth and equivariant under
$A\mapsto RAR^T$.

For $B\in\operatorname{Sym}_0(n)$, differentiation gives
\[
 \begin{aligned}
 q_{D_Ag(s,0)[B]}
 &=P_2D_v\mathscr F(s,0)
       \bigl[q_B+D_Aw(s,0)[B]\bigr]\\
 &=P_2\mathscr L_sq_B=\ell_2(s)q_B,
 \end{aligned}
\]
where we used $D_Aw(s,0)=0$ and the action of
$\mathscr L_s$ on $\mathcal H_2$. Thus
\begin{equation}\label{ordine1}
 D_Ag(s,0)=\ell_2(s)I.
\end{equation}

To compute the quadratic term of $g$, we first expand
the normalized shrinker defect.
 
\begin{lemma}
\label{lem:secondvariation}
For fixed $s\in I$ and $r\ge2$, the normalized defect has the expansion
\[
 \mathscr F(s,v)=\mathscr L_sv+\mathscr B_s^{(2)}(v,v)
                  +O_{Y_{r,s}}\bigl(\|v\|_{X_r}^3\bigr),
 \qquad
 \mathscr B_s^{(2)}=\tfrac12D_v^2\mathscr F(s,0).
\]
The remainder is locally uniform in $s$. The bilinear coefficient is
\begin{align}
 \mathscr B_s^{(2)}(u,v)={}&-\frac{n+s}{4}
 \left(u\frac{s}{\lambda_1(s)}\mathcal I_sv
       +v\frac{s}{\lambda_1(s)}\mathcal I_su\right)
 +\frac{n-2-s}{4}\frac{s}{\lambda_1(s)}\mathcal I_s(uv)\notag\\
 &+\frac{s(1+s)}2uv+\frac12\nabla u\cdot\nabla v.
 \label{eq:Bformula}
\end{align}
\end{lemma}
\begin{proof} 
Recall from \eqref{eq:staticF}
that \[
\mathscr F(s,v)
=
\frac{H_s(E_{1+v},(1+v(\theta))\theta)}{H_s(\mathbb S^{n-1})}
-\frac{1+v}{\sqrt{1+|\nabla v|^2/(1+v)^2}}.
\] For fixed $s$, the second derivative of the normalized shrinker
map with respect to $v$ at $v=0$ is
\begin{align}
 D_v^2\mathscr F(s,0)[u,v]
 ={}&-\frac{n+s}{2H_s(\mathbb S^{n-1})}
 \bigl(u\mathcal I_s v+v\mathcal I_s u\bigr)
 +\frac{n-2-s}{2H_s(\mathbb S^{n-1})}\mathcal I_s(uv)
 \notag\\
 &+s(1+s)uv
 +\nabla u\cdot
  \nabla v.
\end{align}
Indeed, the curvature contribution follows from the second
radial derivative at the unit sphere obtained in \eqref{eq:second-radial-curvature-Rs}, divided by $H_s(\mathbb S^{n-1})$.
The support function is $1+v-\tfrac12|\nabla v|^2+O(v^3)$.
This gives  \eqref{eq:Bformula}, recalling that $\frac{1}{H_s(\mathbb S^{n-1})}=\frac{s}{\lambda_1(s)}$. 
\end{proof} \begin{proposition} \label{lem:quadratic}
For $A\in\operatorname{Sym}_0(n)$, set
\[
 (A^2)_0=A^2-\frac{\operatorname{tr}(A^2)}n I.
\]
Then
\begin{equation}\label{eq:quadraticprojection}
 P_2D_r^2 H_s(1)[q_A,q_A]
 =-\frac{4s(n+s)(n+s+1)}{(n-s)(n+4)}H_s(\mathbb S^{n-1})q_{(A^2)_0}.
\end{equation}
Consequently the coefficient of $q_{(A^2)_0}$ in
$P_2\mathscr B_s^{(2)}(q_A,q_A)$ is
\begin{equation}\label{eq:bnfull}
 b_n(s)=\frac2{n+4}\left[n-\frac{s(n+s)(n+s+1)}{n-s}\right],
\end{equation}
and
\begin{equation}\label{eq:ellipsoid_bn}
 b_n=b_n(\sigma)=-\frac{2(1+\sigma)}{n+4}.
\end{equation}
In particular,
\begin{equation}\label{eq:quadraticbilinearprojection}
 P_2\mathscr B_\sigma^{(2)}(q_A,q_A)
 =-\frac{2(1+\sigma)}{n+4}q_{(A^2)_0}.
\end{equation}
For $n=2$, this projection vanishes identically.
\end{proposition}
\begin{proof}
The spherical moment identities are
\begin{equation}\label{eq:harmonicmoments}
 P_2(q_A^2)=\frac4{n+4}q_{(A^2)_0},\qquad
 P_2(\theta^TA^2\theta)=q_{(A^2)_0}.
\end{equation}
They follow from the rotationally invariant fourth and sixth moments,
or from harmonic decomposition of the degree-four polynomial
$(x^TAx)^2$. Since
$|\nabla q_A|^2=4(\theta^TA^2\theta-q_A^2)$, the support contribution
is $\tfrac12P_2|\nabla q_A|^2=2nq_{(A^2)_0}/(n+4)$.
In \eqref{eq:second-radial-curvature-Rs}, specialized to radius one, harmonic-degree preservation gives
\[
 P_2D^2H_s(1)[q_A,q_A]
 =\left[-\frac{n+3s+2}{2}\lambda_2+s(1+s)H_s(\mathbb S^{n-1})\right]
 \frac4{n+4}q_{(A^2)_0}.
\]
Using
\[
 n(n+3s+2)-(1+s)(n-s)=(n+s)(n+s+1),
\]
substitute $\lambda_2=2nsH_s(\mathbb S^{n-1})/(n-s)$ to obtain
\eqref{eq:quadraticprojection}; adding the support contribution yields
\eqref{eq:bnfull}. At the threshold use
$2n\sigma/(n-\sigma)=1+\sigma$ and
$n=\sigma(1+\sigma)/(1-\sigma)$ to simplify to
\eqref{eq:ellipsoid_bn}; explicitly, the numerator at the threshold is
\[
 -(1+\sigma)(n+\sigma+2)+2n
 =-[\sigma^2+(n+1)\sigma-n]-2(1+\sigma)=-2(1+\sigma).
\]
For $n=2$, $(A^2)_0=0$ identically.
\end{proof}
We can now determine the quadratic expansion of $g$. 
\begin{theorem}\label{thm:quadraticLS}
The jointly smooth reduction satisfies
\begin{equation}\label{eq:matrixnormalform}
 g(s,A)=\ell_2(s)A-\frac{2(1+\sigma)}{n+4}(A^2)_0
 +O\bigl(|s-\sigma||A|^2+|A|^3\bigr).
\end{equation}
For \(n=2\), \((A^2)_0=0\) for every
\(A\in\operatorname{Sym}_0(2)\); hence the quadratic term
vanishes identically.
\end{theorem}
\begin{proof}
By \eqref{eq:reduceddefinition} and Lemma~\ref{lem:secondvariation},
\[
 \begin{aligned}
 q_{g(s,A)}
 ={}&P_2\mathscr L_sq_A+P_2\mathscr L_sw
       +P_2\mathscr B_s^{(2)}(q_A,q_A)\\
    &+2P_2\mathscr B_s^{(2)}(q_A,w)
       +P_2\mathscr B_s^{(2)}(w,w)+O(|A|^3).
 \end{aligned}
\]
Since $P_2w=0$ and $\mathscr L_s$ preserves harmonic degree,
$P_2\mathscr L_sw=\ell_2(s)P_2w=0$.
Moreover, $w=O(|A|^2)$ implies that the two mixed and
quadratic terms involving $w$ are $O(|A|^3)$ and $O(|A|^4)$,
respectively. Proposition~\ref{lem:quadratic} therefore gives
\[
 \begin{aligned}
 q_{g(s,A)}
 &=\ell_2(s)q_A+P_2\mathscr B_s^{(2)}(q_A,q_A)+O(|A|^3)\\
 &=\ell_2(s)q_A+b_n(s)q_{(A^2)_0}+O(|A|^3).
 \end{aligned}
\]
Under the identification $A\mapsto q_A$, this is
\[
 g(s,A)=\ell_2(s)A+b_n(s)(A^2)_0+O(|A|^3),
\]
uniformly for $s$ near $\sigma$.
Since $b_n(s)-b_n(\sigma)=O(|s-\sigma|)$ and
$|(A^2)_0|\le C|A|^2$, replacing $b_n(s)$ by $b_n(\sigma)$
gives the stated remainder.
For $n=2$, the eigenvalues of $A$ are $a,-a$, so
$A^2=a^2I$ and $(A^2)_0=0$.
\end{proof}

In dimension $n=2$ the quadratic term in \eqref{eq:matrixnormalform}
vanishes. In every even dimension $n=2p$, it also vanishes along
the balanced direction, consisting of multiples of the matrix   
\begin{equation}
 A_{\rm b}=\operatorname{diag}(I_p,-I_p),\qquad A_{\rm b}^2=I.
 \label{eq:balancedmatrix}
\end{equation}
We now compute the leading
cubic term in these directions.
 
\begin{proposition} 
\label{thm:balanced}
Let $n=2p\geq2$, $A_{\rm b}$ as in \eqref{eq:balancedmatrix}. 
For $(s,a)$ near $(\sigma,0)$, the reduced map satisfies
\begin{equation}
 g(s,aA_{\rm b})=
 \left[\ell_2(s)a+\gamma_n a^3
 +O\bigl(|s-\sigma||a|^3+|a|^5\bigr)\right]A_{\rm b}.
\label{eq:balancedcubicnormalform}
\end{equation}
Writing
\[
 Q_5 =2\sigma^5-7\sigma^4+12\sigma^3-9\sigma^2+2\sigma+4,
\]
the coefficient is
\begin{equation}
 \gamma_n=
 \frac{2(\sigma-1)(\sigma+1)(3\sigma-1)Q_5}
 {(\sigma-3)(\sigma-2)(\sigma^2-3\sigma+4)(\sigma^2-\sigma+1)}<0.
 \label{eq:cubicclosed}
\end{equation}
In particular, for $n=2$   parametrizing $\mathbb S^1$ by
$(\cos\theta,\sin\theta)$ so that
$ 
Y=\cos2\theta$, we obtain 
\begin{equation}
 \sigma=\frac{\sqrt{17}-3}{2},\qquad
 \gamma_2=\frac{3585-889\sqrt{17}}{156}<0.
 \label{eq:planarcubiccoefficient}
\end{equation}
\end{proposition}

 \begin{proof}
We first use symmetry to identify the form of the reduced
equation. We then compute its cubic coefficient by comparing
the range-corrected graph with an explicit family of ellipsoids.

{\bf Step 1: proof of \eqref{eq:balancedcubicnormalform}.}

We first observe that the matrices in $\operatorname{Sym}_0(n)$
fixed under conjugation by $O(p)\times O(p)$ are precisely
the multiples of $A_{\rm b}$.
Write
\[
 A=\begin{pmatrix}C&D\\D^T&E\end{pmatrix},
 \qquad
 R=\begin{pmatrix}U&0\\0&V\end{pmatrix},
 \qquad U,V\in O(p),
\]
where $O(p)=\{U\in\mathbb R^{p\times p}:U^TU=I_p\}$,
and assume that $RAR^T=A$ for every such $R$.
Taking $U=-I_p$ and $V=I_p$ gives $D=-D$, hence $D=0$.
The remaining conditions are
\[
 UCU^T=C,\qquad VEV^T=E
 \qquad\text{for all }U,V\in O(p).
\]
Invariance under coordinate reflections forces $C$ and $E$
to be diagonal, while invariance under coordinate permutations
forces their diagonal entries to be equal. Thus
$C=\alpha I_p$ and $E=\beta I_p$.
Since $\operatorname{tr}A=p(\alpha+\beta)=0$, we have
$\beta=-\alpha$, and therefore $A=\alpha A_{\rm b}$.
Conversely, every multiple of $A_{\rm b}$ is fixed by these
conjugations.

Since $aA_{\rm b}$ is fixed by the subgroup, equivariance gives
\[
 Rg(s,aA_{\rm b})R^T
 =g(s,RaA_{\rm b}R^T)
 =g(s,aA_{\rm b})
 \qquad\text{for every }R\in O(p)\times O(p).
\]
Hence $g(s,aA_{\rm b})\in\mathbb R A_{\rm b}$, so the reduced
equation along this line is scalar.
Moreover, the orthogonal transformation interchanging the two
blocks sends $A_{\rm b}$ to $-A_{\rm b}$. Applying equivariance
once more gives
\[
 g(s,-aA_{\rm b})=-g(s,aA_{\rm b}).
\]
The scalar coefficient of $A_{\rm b}$ is therefore odd in $a$.

By joint smoothness and $D_Ag(s,0)=\ell_2(s)I$, its Taylor
expansion consists of the linear term $\ell_2(s)a$, a cubic
term, and a remainder $O(|a|^5)$, uniformly for $s$ near
$\sigma$. The cubic coefficient depends smoothly on $s$;
denote its value at $s=\sigma$ by $\gamma_n$.
Replacing that coefficient by $\gamma_n$ produces an error
$O(|s-\sigma||a|^3)$.
This proves the form of \eqref{eq:balancedcubicnormalform}.
It remains to compute $\gamma_n$.

 {\bf Step 2: computation of $\gamma_n$.}

Recall that $Y=q_{A_{\rm b}}$. The spherical moments give
\[
 m_2=\langle Y^2\rangle=\frac{2}{n+2},\qquad
 \langle Y^4\rangle=\frac{12}{(n+2)(n+6)},\qquad
 c_3=\frac{\langle Y^4\rangle}{\langle Y^2\rangle}
     =\frac{6}{n+6}.
\]
Since $A_{\rm b}^2=I$, the quadratic projection formula yields
$P_2(Y^2)=0$; hence $Y_4:=Y^2-m_2\in\mathcal H_4$.
By the fixed-subspace argument in Step~1, $P_2(Y^3)$ is a
multiple of $Y$, and taking its inner product with $Y$
determines the coefficient. Thus
\[
 P_2(Y^3)=c_3Y,\qquad
 P_2(YY_4)=(c_3-m_2)Y,\qquad
 |\nabla Y|^2=4(1-Y^2).
\]

For $0<s<1$, define
\begin{equation}\label{eq:cosinetransform}
 (\mathcal T_s f)(\nu)
 =\frac{\displaystyle\int_{\Sph}
          |\nu\cdot\omega|^{-s}f(\omega)\dd\omega}
        {\displaystyle\int_{\Sph}
          |\nu\cdot\omega|^{-s}\dd\omega}.
\end{equation}
This rotationally invariant operator acts on $\mathcal H_2$
and $\mathcal H_4$ with multipliers
\begin{equation}\label{eq:cosineeigen}
 \tau_2=-\frac{s}{n-s},\qquad
 \tau_4=\frac{s(s+2)}{(n-s)(n-s+2)}.
\end{equation}
These follow by evaluating zonal harmonics at $\nu=e_1$, using
\[
 \frac{\int_{\Sph}|\omega_1|^{-s}\omega_1^{2j}\dd\omega}
      {\int_{\Sph}|\omega_1|^{-s}\dd\omega}
 =\frac{((1-s)/2)_j}{((n-s)/2)_j},\qquad j=1,2.
\]
Consequently,
\[
 \mathcal T_sY=\tau_2Y,\qquad
 \mathcal T_sY^2=\tau_4Y^2+(1-\tau_4)m_2,\qquad
 P_2\mathcal T_sY^3=\tau_2c_3Y.
\]

Consider the ellipsoids
\[
 Q_a=I+2aA_{\rm b},\qquad
 \mathbb E_a=\{x\in\R^n:x^TQ_a^{-1}x<1\},\qquad |a|<\tfrac12.
\]
Their support function and the boundary point with outer normal
$\nu$ are
\[
 h_a(\nu)=\sqrt{1+2aY(\nu)},\qquad
 x(\nu)=\frac{Q_a\nu}{h_a(\nu)}.
\]
The inward chord in direction $\omega\cdot\nu<0$ has length
\[
 L(\omega)=\frac{2|\nu\cdot\omega|}
                 {h_a(\nu)\,\omega^TQ_a^{-1}\omega}.
\]
The tangent half-space $\Pi$ contains $\mathbb E_a$ and has
zero curvature at $x=x(\nu)$. Subtracting its contribution
and using polar coordinates centered at $x$ gives
\[
 \begin{aligned}
 H_s(\mathbb E_a,x)
=2\int_{\Pi\setminus\mathbb E_a}|x-y|^{-n-s}\dd y =2\int_{\{\omega\cdot\nu<0\}}
\int_{L(\omega)}^\infty t^{-1-s}\dd t\dd\omega =\frac2s\int_{\{\omega\cdot\nu<0\}}L(\omega)^{-s}\dd\omega.
 \end{aligned}
\]
The chord formula extends evenly in $\omega$, and its inverse
$s$th power is integrable because $s<1$. Comparing with the
unit ball and using $Q_a^{-1}=(I-2aA_{\rm b})/(1-4a^2)$ yields
\begin{equation}\label{eq:ellipsoididentity}
 \frac{H_s(\mathbb E_a,x(\nu))}{H_s(\mathbb S^{n-1})}
 =(1+2aY(\nu))^{s/2}(1-4a^2)^{-s}
   \mathcal T_s[(1-2aY)^s](\nu).
\end{equation}

From now on, $s=\sigma$ and all coefficients are evaluated at
this value. Expanding \eqref{eq:ellipsoididentity} and subtracting
$h_a$ gives a vanishing linear term, since $\ell_2(\sigma)=0$.
The quadratic term is $\overline f+f_4Y_4$, and the degree-two
projection of the cubic term is $e_3Y$. Using the identities
above and $n=\sigma(1+\sigma)/(1-\sigma)$, their coefficients
simplify to
\begin{equation}\label{eq:ellipsoidsecond}
 f_4=-\frac{(\sigma-1)(3\sigma-1)}{2(\sigma^2-\sigma+1)},
 \qquad
 \overline f=\frac{(\sigma+1)^3}{\sigma^2-\sigma+2},
\end{equation}
and
\begin{equation}\label{eq:ellipsoidthird}
 e_3=\frac{(\sigma+1)
       (\sigma^4-9\sigma^3+16\sigma^2-9\sigma+5)}
      {(\sigma-3)(\sigma-2)(\sigma^2-\sigma+1)}.
\end{equation}

We next express this defect in radial coordinates. The radial
function and the change from the normal direction $\nu$ to
$\theta=x/|x|$ satisfy
\begin{equation}\label{eq:ellipsoidradial}
 \begin{aligned}
 r_a(\theta)
 =\sqrt{1-4a^2}\,(1-2aY(\theta))^{-1/2}=1+aY(\theta)+a^2\bigl(\tfrac32Y(\theta)^2-2\bigr)+O(a^3),
 \end{aligned}
\end{equation}
\[
 \theta=\frac{Q_a\nu}{|Q_a\nu|}
       =\nu+a\nabla Y(\nu)+O(a^2),\qquad
 \nu=\theta-a\nabla Y(\theta)+O(a^2).
\]
The defect is unchanged at corresponding boundary points:
\[
 \mathscr F(\sigma,r_a-1)(\theta(\nu))
 =\frac{H_\sigma(\mathbb E_a,x(\nu))}
        {H_\sigma(\mathbb S^{n-1})}-h_a(\nu).
\]
Since its linear term vanishes, the only cubic correction
comes from composing the quadratic term with $\nu(\theta)$.
As $\nabla Y_4=2Y\nabla Y$, this correction is
$-2a^3f_4Y|\nabla Y|^2$. Therefore
\[
 \begin{aligned}
 P_2(Y|\nabla Y|^2)&=4(1-c_3)Y,\\
 P_2\mathscr F(\sigma,r_a-1)
 &=a^3\bigl[e_3-8f_4(1-c_3)\bigr]Y+O(a^4).
 \end{aligned}
\]

It remains to replace the ellipsoidal quadratic correction
$\tfrac32Y^2-2$ by the correction determined by the range
equation. Write
\[
 w(\sigma,aA_{\rm b})=a^2w_2+O(a^3),\qquad P_2w_2=0.
\]
Since $P_2\mathscr B_\sigma^{(2)}(Y,Y)=0$, comparison of
the quadratic terms gives
\[
 \begin{aligned}
 \mathscr L_\sigma w_2+\mathscr B_\sigma^{(2)}(Y,Y)&=0,\\
 \mathscr L_\sigma\bigl(\tfrac32Y^2-2\bigr)
       +\mathscr B_\sigma^{(2)}(Y,Y)&=\overline f+f_4Y_4.
 \end{aligned}
\]
Let
\begin{equation}\label{eq:t4}
 t_4=\sigma\sum_{j=0}^{3}
       \frac{((n+\sigma)/2)_j}{((n-\sigma)/2)_j},
 \qquad \ell_4=t_4-(1+\sigma)>0.
\end{equation}
Subtracting the two equations and using
$\mathscr L_\sigma1=-(1+\sigma)$ and
$\mathscr L_\sigma Y_4=\ell_4Y_4$ yields
\begin{equation}\label{eq:rangecorrectionbalanced}
 w_2-\bigl(\tfrac32Y^2-2\bigr)
 =\frac{\overline f}{1+\sigma}-\frac{f_4}{\ell_4}Y_4.
\end{equation}
Here the inversion is on the complement of $\mathcal H_2$,
which contains both quadratic corrections.

The bilinear formula \eqref{eq:Bformula} gives
\begin{align}
 P_2\mathscr B_\sigma^{(2)}(Y,1)
 &=-\frac{1+\sigma}{2}Y,\notag\\
 P_2\mathscr B_\sigma^{(2)}(Y,Y_4)
 &=B_{24}Y,\qquad
 B_{24}=\left[n+2-\frac{1+\sigma}{2}
                  -\frac{n+\sigma}{4}t_4\right](c_3-m_2).
 \label{eq:crossB}
\end{align}
For these identities, use
$\sigma\lambda_2(\sigma)/\lambda_1(\sigma)=1+\sigma$ and
$P_2(\nabla Y\cdot\nabla Y_4)=2(n+2)P_2(YY_4)$;
the latter follows from the spherical Laplacian product rule.

Both radial graphs have first-order term $aY$, so their
third-derivative contributions to the cubic coefficient
coincide. Their linear contributions vanish after projection
because $P_2\mathscr L_\sigma=0$. Thus their projected cubic
coefficients differ only by
\[
 2P_2\mathscr B_\sigma^{(2)}
       \bigl(Y,w_2-(\tfrac32Y^2-2)\bigr)
 =\left[-\overline f-\frac{2f_4}{\ell_4}B_{24}\right]Y.
\]
Combining this range correction with the coordinate correction
gives
\begin{equation}\label{eq:cubicunsimplified}
 \gamma_n=e_3-8f_4(1-c_3)-\overline f
                 -\frac{2f_4}{\ell_4}B_{24}.
\end{equation}
No third-order radial corrections are needed, since their
images under $\mathscr L_\sigma$ have zero degree-two projection.

Finally, substituting $n=\sigma(1+\sigma)/(1-\sigma)$ into
\eqref{eq:t4} gives
\[
 \ell_4=\frac{\sigma^2-3\sigma+4}
 { (\sigma^2-2\sigma+2)(\sigma^2-\sigma+1)}.
\]
Together with \eqref{eq:ellipsoidsecond}--\eqref{eq:ellipsoidthird}
and \eqref{eq:crossB}, this reduces
\eqref{eq:cubicunsimplified} to \eqref{eq:cubicclosed}.
Its sign follows from $1/3<\sigma<1$ and the identity
\[ 
 Q_5=4(1-\sigma)^5+22\sigma(1-\sigma)^4
        +39\sigma^2(1-\sigma)^3 +37\sigma^3(1-\sigma)^2
        +18\sigma^4(1-\sigma)+4\sigma^5>0.
\]
The denominator in \eqref{eq:cubicclosed} is positive, while
its numerator is negative; hence $\gamma_n<0$.

{\bf Step 3: case $n=2$.}

Finally, for $n=2$,
\[
 m_2=\frac12,\qquad c_3=\frac34,\qquad
 Y_4=\frac12\cos4\theta.
\]
The coordinate correction is then $2f_4$.
Writing $q=\sqrt{17}$, the quantities entering
\eqref{eq:cubicunsimplified} become
\[
 \ell_4=\frac{39+9q}{26},
 \qquad
 f_4=\frac{-35+11q}{52},
 \qquad
 \overline f=\frac{5+3q}{8},
\]
\[
 e_3=\frac{-127+135q}{208},
 \qquad
 B_{24}=\frac{121-25q}{208}.
\]
Thus
\[
 \gamma_2
 =e_3-2f_4-\overline f-\frac{2f_4}{\ell_4}B_{24}
 =\frac{3585-889\sqrt{17}}{156}<0.
\]
\end{proof}

\section{Nearly spherical shrinkers}\label{sec:branches}
We now construct nonspherical shrinkers bifurcating from
the sphere at the first threshold. Unbalanced branches occur
in dimensions $n\ge3$, whereas balanced branches occur in
every even dimension.

\subsection{Unbalanced block symmetries: transverse branches}
Let $n\geq 3$.

For $1\leq p<n/2$ put $q=n-p$, write $\theta=(\theta',\theta'')\in\mathbb R^p\times\mathbb R^q$, and define
\begin{equation}
 A_p=\operatorname{diag}\left(\frac qn I_p,-\frac pn I_q\right),
 \qquad Y_p(\theta)=|\theta'|^2-\frac pn,
 \qquad (A_p^2)_0=\frac{n-2p}{n}A_p.
\label{eq:blockmatrix}
\end{equation}
The fixed subspace of $\operatorname{Sym}_0(n)$ under $O(p)\times O(q)$ is exactly $\mathbb R A_p$.

\begin{theorem} \label{thm:unbalanced}
For each $1\leq p<n/2$ there is a smooth branch of centrally
symmetric, strictly convex shrinkers, nonspherical for $t\ne0$,
\begin{equation}
 \partial\mathcal E_{p,t}
 =\left\{R_{s_p(t)}\bigl[1+tY_p(\theta)
                 +w(s_p(t),tA_p)(\theta)\bigr]\theta:\theta\in\Sph\right\},
 \label{eq:unbalancedbranch}
\end{equation}
for $|t|$ sufficiently small. Here
\begin{equation}
 s_p(0)=\sigma,\qquad
 s_p'(0)=\frac{2(1+\sigma)(n-2p)}{n(n+4)d_n}>0.
 \label{eq:unbalancedslope}
\end{equation}
Every nonzero member satisfies $H_{s_p(t)}=x\cdot\nu/(1+s_p(t))$ and is not a sphere. The branch crosses the threshold: positive $t$ lies on the $s>\sigma$ side and negative $t$ on the $s<\sigma$ side.
\end{theorem}
\begin{proof}
Equivariance implies $g(s,tA_p)=g_p(s,t)A_p$. Since $g_p(s,0)=0$, its quotient by $t$ extends smoothly:
\[
 \frac{g_p(s,t)}t=\ell_2(s)+b_n\frac{n-2p}{n}t              +O(t^2+|s-\sigma||t|).
\]
 Since the derivative with respect to $s$ in $(\sigma,0)$
is $d_n>0$, the scalar implicit function theorem produces
a unique smooth function $s_p(t)$ solving the divided
equation, with $s_p(0)=\sigma$. Implicit differentiation
gives \eqref{eq:unbalancedslope}.

The corresponding physical radial height is
$u_{p,t}=R_{s_p(t)}\,[tY_p+w(s_p(t),tA_p)]$.
The range solution is smooth, even, and invariant under the same block group. Since $P_2w(s_p(t),tA_p)=0$, the radial perturbation satisfies
\[
 P_2\bigl[tY_p+w(s_p(t),tA_p)\bigr]=tY_p\ne0
 \qquad\text{for }t\ne0.
\]
A sphere centered at the origin has a constant radial
function and therefore zero degree-two projection.
Moreover, the constructed hypersurface is centrally
symmetric about the origin, so any spherical member
would necessarily be centered there.
Thus every member with $t\ne0$ is nonspherical. Equation \eqref{eq:staticF} gives the shrinker equation after multiplication by $R_s$. Uniform $C^2$ smallness implies positive classical principal curvatures and embeddedness. 
\end{proof}

Since $s_p'(0)>0$, the inverse function theorem gives a smooth
inverse $t_p(s)$ for $|s-\sigma|<\varepsilon_p$, where
$\varepsilon_p>0$ and $t_p(\sigma)=0$. Taylor expansion and the estimate
$w(s,A)=O(|A|^2)$ from Lemma~\ref{lem:range} yield
\begin{equation}\label{eq:subcriticalbranchrate}
 t_p(s)=\frac{s-\sigma}{s_p'(0)}+O((s-\sigma)^2),\qquad
 \|t_p(s)Y_p+w(s,t_p(s)A_p)\|_{C^k}
 \le C_k|s-\sigma|
\end{equation}
for every fixed $k$. Thus, for every $s\ne\sigma$ sufficiently
close to $\sigma$, the profile $\mathcal E_{p,t_p(s)}$ is a smooth,
embedded, strictly convex nonspherical shrinker. These profiles
approach the round sphere as $s\to\sigma$ from either side.

\subsection{Balanced branches}

For even $n=2p\geq 2$, use $A_{\rm b}$ and $Y$ from
\eqref{eq:balancedmatrix}.  Proposition~\ref{thm:balanced} supplies
the cubic reduced equation \eqref{eq:balancedcubicnormalform}
and its negative coefficient. We now use this equation to construct
the balanced branch.

\begin{theorem}\label{thm:balancedexist}
Let $n=2p$. There is a smooth, even function $s_{\rm b}(a)$
defined near zero such that
\[
 s_{\rm b}(a)=\sigma-\frac{\gamma_n}{d_n}a^2+O(a^4),
 \qquad s_{\rm b}(a)>\sigma\quad\text{for }a\ne0.
\]
The hypersurfaces
\[
 \left\{R_{s_{\rm b}(a)}
 \bigl[1+aY(\theta)+w(s_{\rm b}(a),aA_{\rm b})(\theta)\bigr]\theta:
 \theta\in\mathbb S^{n-1}\right\}
\]
are smooth, strictly convex, $O(p)\times O(p)$-invariant
shrinkers for all sufficiently small $a$.
They satisfy $H_{s_{\rm b}(a)}=x\cdot\nu/(1+s_{\rm b}(a))$
and are nonspherical whenever $a\ne0$.
\end{theorem}
\begin{proof}
By equivariance, $g(s,aA_{\rm b})$ is a scalar multiple
of $A_{\rm b}$, and this scalar is odd in $a$.
Its quotient by $a$ extends smoothly and evenly to $a=0$.
Proposition~\ref{thm:balanced} gives the divided equation
\[
 \ell_2(s)+\gamma_n a^2
 +O\bigl(|s-\sigma|a^2+a^4\bigr)=0.
\]
Its $s$-derivative at $(\sigma,0)$ is $d_n>0$.
The implicit function theorem therefore gives a unique
smooth function $s=s_{\rm b}(a)$. Evenness and uniqueness
imply $s_{\rm b}(-a)=s_{\rm b}(a)$, and the displayed
equation yields the asserted expansion.
The sign follows from $\gamma_n<0$ and $d_n>0$.

The range solution is smooth and $O(a^2)$ in every $C^k$ norm,
so the resulting hypersurfaces are strictly convex for small $a$.
Multiplication by $R_{s_{\rm b}(a)}$ converts the dimensionless
equation into the stated physical shrinker equation.
Finally, $P_2[aY+w(s_{\rm b}(a),aA_{\rm b})]=aY\ne0$
when $a\ne0$. Central symmetry forces any spherical member
to be centered at the origin, where its radial function would
be constant. Thus the nonzero members are nonspherical.
\end{proof}

The same quantities determine the second-order shape of the balanced profiles.

\begin{proposition}
\label{prop:balancedprofile}
Let $n=2p$ be even, $\sigma=s_*(n)$, and use the amplitude convention \eqref{eq:balancedmatrix}. Then $w(\sigma,aA_{\rm b})=a^2w_2+O(a^3)$ in every $X_r$, where
\begin{equation}
 w_2=\frac{\sigma(1-\sigma)}{\sigma^2-\sigma+2}
 +\frac{3\sigma^4-10\sigma^3+18\sigma^2-19\sigma+14}{2(\sigma^2-3\sigma+4)}\,Y_4,
 \qquad Y_4=Y^2-\frac{2}{n+2}\in\mathcal H_4.
 \label{eq:balancedw2}
\end{equation}
Consequently, along the balanced branch $s=s_{\rm b}(a)$ of Theorem~\ref{thm:balancedexist}, the radial profile is
$1+aY+a^2w_2+O(a^3)$ in every $C^k$ norm. In particular, for $a\ne0$ sufficiently small the balanced shrinkers are not ellipsoids.
\end{proposition}
\begin{proof}
By \eqref{eq:ellipsoidradial}, $ \tfrac32Y^2-2=\tfrac32Y_4+\tfrac32m_2-2$. Inserting this into \eqref{eq:rangecorrectionbalanced} gives
\[
 w_2=\Bigl[\tfrac32m_2-2+\frac{\overline f}{1+\sigma}\Bigr]
 +\Bigl[\tfrac32-\frac{f_4}{\ell_4}\Bigr]Y_4.
\]
At the threshold $n=\sigma(1+\sigma)/(1-\sigma)$, so that $m_2=2(1-\sigma)/(\sigma^2-\sigma+2)$. With the values of $\ell_4$, $f_4$ and $\overline f$ listed at the end of the proof of Proposition~\ref{thm:balanced}, both brackets simplify to the coefficients in \eqref{eq:balancedw2}. The expansion along the branch follows because
$s_{\rm b}(a)-\sigma=O(a^2)$ and
$w(s,aA_{\rm b})=a^2w_2(s)+O(a^3)$, with $w_2(s)$ smooth
in $s$ by Lemma~\ref{lem:range}.

An ellipsoid invariant under $O(p)\times O(p)$ and centered at the origin has radial function $C(1+tY)^{-1/2}$ for some constants $C,t$. If such an ellipsoid has the same radial profile as a member
of the branch, comparison of its constant and degree-two
components gives $C=1+O(a^2)$ and $t=-2a+O(a^3)$.
Expanding $C(1+tY)^{-1/2}$ then shows that its degree-four
component is $\tfrac32a^2Y_4+O(a^3)$. For the shrinkers this coefficient is $\tfrac32-f_4/\ell_4$, and $f_4=-(\sigma-1)(3\sigma-1)/(2(\sigma^2-\sigma+1))\neq0$ because $\sigma>1/3$. Thus the degree-four components differ by
$-(f_4/\ell_4)a^2Y_4+O(a^3)$, which cannot vanish for
sufficiently small nonzero $a$.
Central symmetry and block invariance force any ellipsoid
coinciding with a branch member to have the form considered above.
This proves the last assertion.
\end{proof}

For $n=2$, parametrizing $\mathbb S^1$ by
$(\cos\theta,\sin\theta)$ gives
\[
Y=\cos2\theta,\qquad
m_2=\frac12,\qquad
Y_4=Y^2-m_2
=\cos^2(2\theta)-\frac12
=\frac12\cos4\theta,
\] and \eqref{eq:balancedw2} reduces to the planar profile in Corollary~\ref{cor:planar}.

\begin{corollary} 
\label{cor:planar}
Let $n=2$ and use the amplitude $a\cos2\theta$, that is, $A=aA_{\rm b}$ with $A_{\rm b}=\operatorname{diag}(1,-1)$. Then
\begin{gather}
 \sigma=\frac{\sqrt{17}-3}{2},\qquad
d_2=\ell_2'(\sigma)=\frac{17+7\sqrt{17}}{16},\qquad \gamma_2=\frac{3585-889\sqrt{17}}{156}
\label{eq:planarexact}
\end{gather}
and the exponent along the branch satisfies
\begin{equation}
 s(a)=\sigma-\frac{\gamma_2}{d_2}\,a^2+O(a^4).\label{eq:planarbranch}
\end{equation}
The range correction is
\begin{equation}
 w(\sigma,aA_{\rm b})=a^2\left[\frac{\sqrt{17}-3}{8}+\frac{145-31\sqrt{17}}{24}\cos4\theta\right]+O(a^3).
 \label{eq:planarw2}
\end{equation}
Up to rotations, this is the only local family of noncircular shrinkers near the configuration $(\sigma,B_1)$. No noncircular shrinker lies sufficiently close to that point with $s\le\sigma$. The assertion is local in both exponent and shape.
\end{corollary}

\begin{proof}
Every trace-free symmetric $2\times2$ matrix is orthogonally conjugate to $a\operatorname{diag}(1,-1)$. Equivariance and uniqueness of the range solution reduce all possible zeros to the balanced scalar equation
\[
 a\left[\ell_2(s)+\gamma_2a^2+O(|s-\sigma|a^2+a^4)\right]=0.
\]
For $a\ne0$, divide by $a$. Since $d_2=8/(2-\sigma)^2-1>0$, the implicit function theorem gives a unique smooth solution $s=s(a)$ near $(\sigma,0)$; block exchange sends $a$ to $-a$, so $s(a)$ is even, and $s(a)=\sigma-(\gamma_2/d_2)a^2+O(a^4)$.

The value of $\gamma_2$ is given by
Proposition~\ref{thm:balanced}, and the displayed value of $d_2$
follows by differentiating $\ell_2(s)$.   Finally \eqref{eq:planarw2} is Proposition~\ref{prop:balancedprofile} at $n=2$.

Since $-\frac{\gamma_2}{d_2}>0$, the branch lies above the threshold. Uniqueness of the range solution and of the scalar implicit branch excludes every other sufficiently small nonzero solution at or below the threshold.
\end{proof}




\begin{thebibliography}{99}
\small

\bibitem{AbreschLanger}
U. Abresch and J. Langer,
\emph{The normalized curve shortening flow and homothetic solutions},
J. Differential Geom. \textbf{23} (1986), no.~2, 175--196.

\bibitem{CFSMW}
X. Cabr\'e, M. M. Fall, J. Sol\`a-Morales, and T. Weth,
\emph{Curves and surfaces with constant nonlocal mean curvature: meeting Alexandrov and Delaunay},
J. Reine Angew. Math. \textbf{745} (2018), 253--280.
\href{https://doi.org/10.1515/crelle-2015-0117}{doi:10.1515/crelle-2015-0117}.

\bibitem{CFWdelaunay}
X. Cabr\'e, M. M. Fall, and T. Weth,
\emph{Delaunay hypersurfaces with constant nonlocal mean curvature},
J. Math. Pures Appl. \textbf{110} (2018), 32--70.
\href{https://doi.org/10.1016/j.matpur.2017.07.005}{doi:10.1016/j.matpur.2017.07.005}.

\bibitem{CFWlattices}
X. Cabr\'e, M. M. Fall, and T. Weth,
\emph{Near-sphere lattices with constant nonlocal mean curvature},
Math. Ann. \textbf{370} (2018), no.~3--4, 1513--1569.
\href{https://doi.org/10.1007/s00208-017-1559-6}{doi:10.1007/s00208-017-1559-6}.

\bibitem{CDNP}
A. Cesaroni, L. De Luca, M. Novaga, and M. Ponsiglione,
\emph{Stability results for nonlocal geometric evolutions and limit cases for fractional mean curvature flows},
Commun. Partial Differential Equations \textbf{46} (2021), no.~7, 1344--1371.
\href{https://doi.org/10.1080/03605302.2021.1875485}{doi:10.1080/03605302.2021.1875485}.

\bibitem{CNshrink}
A. Cesaroni and M. Novaga,
\emph{Symmetric self-shrinkers for the fractional mean curvature flow},
J. Geom. Anal. \textbf{30} (2020), no.~4, 3698--3715.
\href{https://doi.org/10.1007/s12220-019-00214-2}{doi:10.1007/s12220-019-00214-2}.

\bibitem{CNgraphs}
A. Cesaroni and M. Novaga,
\emph{Fractional mean curvature flow of Lipschitz graphs},
Manuscripta Math. \textbf{170} (2023), 427--451.
\href{https://doi.org/10.1007/s00229-022-01371-5}{doi:10.1007/s00229-022-01371-5}.

\bibitem{CNball}
A. Cesaroni and M. Novaga,
\emph{Stability of the ball under volume preserving fractional mean curvature flow},
Adv. Calc. Var. \textbf{17} (2024), no.~2, 503--520.
\href{https://doi.org/10.1515/acv-2022-0027}{doi:10.1515/acv-2022-0027}.

\bibitem{CMP}
A. Chambolle, M. Morini, and M. Ponsiglione,
\emph{Nonlocal curvature flows},
Arch. Ration. Mech. Anal. \textbf{218} (2015), 1263--1329.
\href{https://doi.org/10.1007/s00205-015-0880-z}{doi:10.1007/s00205-015-0880-z}.

\bibitem{CNR}
A. Chambolle, M. Novaga, and B. Ruffini,
\emph{Some results on anisotropic fractional mean curvature flows},
Interfaces Free Bound. \textbf{19} (2017), no.~3, 393--415.
\href{https://doi.org/10.4171/IFB/387}{doi:10.4171/IFB/387}.

\bibitem{CSVneck}
E. Cinti, C. Sinestrari, and E. Valdinoci,
\emph{Neckpinch singularities in fractional mean curvature flows},
Proc. Amer. Math. Soc. \textbf{146} (2018), no.~6, 2637--2646.
\href{https://doi.org/10.1090/proc/14002}{doi:10.1090/proc/14002}.

\bibitem{CSVvolume}
E. Cinti, C. Sinestrari, and E. Valdinoci,
\emph{Convex sets evolving by volume-preserving fractional mean curvature flows},
Anal. PDE \textbf{13} (2020), no.~7, 2149--2171.
\href{https://doi.org/10.2140/apde.2020.13.2149}{doi:10.2140/apde.2020.13.2149}.

\bibitem{CFMN}
G. Ciraolo, A. Figalli, F. Maggi, and M. Novaga,
\emph{Rigidity and sharp stability estimates for hypersurfaces with constant and almost-constant nonlocal mean curvature},
J. Reine Angew. Math. \textbf{741} (2018), 275--294.
\href{https://doi.org/10.1515/crelle-2015-0088}{doi:10.1515/crelle-2015-0088}.

\bibitem{CR}
M. G. Crandall and P. H. Rabinowitz,
\emph{Bifurcation from simple eigenvalues},
J. Funct. Anal. \textbf{8} (1971), no.~2, 321--340.
\href{https://doi.org/10.1016/0022-1236(71)90015-2}{doi:10.1016/0022-1236(71)90015-2}.

\bibitem{FFMMM}
A. Figalli, N. Fusco, F. Maggi, V. Millot, and M. Morini,
\emph{Isoperimetry and stability properties of balls with respect to nonlocal energies},
Commun. Math. Phys. \textbf{336} (2015), no.~1, 441--507.
\href{https://doi.org/10.1007/s00220-014-2244-1}{doi:10.1007/s00220-014-2244-1}.
For the corrected version used for the second-variation formulas, see arXiv:1403.0516v2.

\bibitem{GageHamilton}
M. Gage and R. S. Hamilton,
\emph{The heat equation shrinking convex plane curves},
J. Differential Geom. \textbf{23} (1986), no.~1, 69--96.
\href{https://doi.org/10.4310/jdg/1214439902}{doi:10.4310/jdg/1214439902}.

\bibitem{Huisken}
G. Huisken,
\emph{Asymptotic behavior for singularities of the mean curvature flow},
J. Differential Geom. \textbf{31} (1990), no.~1, 285--299.
\href{https://doi.org/10.4310/jdg/1214444099}{doi:10.4310/jdg/1214444099}.

\bibitem{Imbert}
C. Imbert,
\emph{Level set approach for fractional mean curvature flows},
Interfaces Free Bound. \textbf{11} (2009), no.~1, 153--176.
\href{https://doi.org/10.4171/IFB/207}{doi:10.4171/IFB/207}.

\bibitem{JLMshort}
V. Julin and D. A. La Manna,
\emph{Short time existence of the classical solution to the fractional mean curvature flow},
Ann. Inst. H. Poincar\'e Anal. Non Lin\'eaire \textbf{37} (2020), no.~4, 983--1016.
\href{https://doi.org/10.1016/j.anihpc.2020.02.007}{doi:10.1016/j.anihpc.2020.02.007}.

\bibitem{JLMconvex}
V. Julin and D. A. La Manna,
\emph{Convergence of the volume preserving fractional mean curvature flow for convex sets},
Ann. Sc. Norm. Super. Pisa Cl. Sci. (5), published online 29 May 2025, 30 pp.
\href{https://doi.org/10.2422/2036-2145.202311_019}{doi:10.2422/2036-2145.202311\_019}.

\bibitem{MWstability}
B.-V. Matioc and C. Walker,
\emph{On the principle of linearized stability in interpolation spaces for quasilinear evolution equations},
Monatsh. Math. \textbf{191} (2020), no.~3, 615--634.
\href{https://doi.org/10.1007/s00605-019-01352-z}{doi:10.1007/s00605-019-01352-z}.

\bibitem{MWgraphs}
B.-V. Matioc and C. Walker,
\emph{The nonlocal fractional mean curvature flow of periodic graphs},
Ann. Sc. Norm. Super. Pisa Cl. Sci. (5) \textbf{26} (2025), no.~1, 91--130.
\href{https://doi.org/10.2422/2036-2145.202208_019}{doi:10.2422/2036-2145.202208\_019}.

\bibitem{NPW}
F. J. Narcowich, P. Petrushev, and J. D. Ward,
\emph{Localized tight frames on spheres},
SIAM J. Math. Anal. \textbf{38} (2006), no.~2, 574--594.
\href{https://doi.org/10.1137/040614359}{doi:10.1137/040614359}.

\bibitem{SV}
M. S\'aez and E. Valdinoci,
\emph{On the evolution by fractional mean curvature},
Commun. Anal. Geom. \textbf{27} (2019), no.~1, 211--249.
\href{https://doi.org/10.4310/CAG.2019.v27.n1.a6}{doi:10.4310/CAG.2019.v27.n1.a6}.

\end{thebibliography}
\end{document}